\documentclass{article}

\usepackage{amsmath,amsfonts,amssymb,amstext,amsthm,amscd,mathrsfs}
\usepackage{epsfig}
\usepackage[all]{xy}

\def\C{{\bf \mathbb{C}}}
\def\N{{\bf \mathbb{N}}}

\def\R{{\bf \mathbb{R}}}

\def\P{{\bf \mathbb{P}}}

\def\m{{\mathfrak{m}}}
\def\M{{\mathfrak{M}}}
\def\X{{\mathfrak{X}}}

\def\z{{\underline{z}}}
\def\a{{\underline{a}}}

\def\elem(#1,#2){  \{ {{#1}\over \overline {\ #2\ }}\}  }

\def\X{{\mathfrak{X}}}

\def\z{{\underline{z}}}
\def\a{{\underline{a}}}

\def \Na{{\mathcal{N}}}

\newtheorem{definition}{Definition}[section]

\newtheorem{remark}[definition]{Remark}
\newtheorem{ex}{Example}[section]

\newenvironment{example}{\begin{ex}\em}{\end{ex}}
\newtheorem{Theorem}[definition]{Theorem}
\newtheorem{corollary}[definition]{Corollary}

\newtheorem{proposition}[definition]{Proposition}
\newtheorem{lema}[definition]{Lemma}
\begin{document}
\title{Specialization to the Tangent Cone and Whitney Equisingularity.}
\author{Arturo Giles Flores}
\maketitle
\pagestyle{myheadings}
\markboth{\rm Arturo Giles Flores}{\rm Specialization to the tangent cone.}
\vskip.3pt 	

\begin{abstract}
     Let $(X,0)$ be a reduced, equidimensional germ of analytic singularity with reduced tangent
 cone $(C_{X,0},0)$. We prove that the absence of exceptional cones is a necessary and sufficient condition 
 for the smooth part $\X^0$ of the specialization to the tangent cone $\varphi: \X \to \C$ 
 to satisfy Whitney's conditions along the parameter axis $Y$. This result is a first step in ge- neralizing
 to higher dimensions L\^e and Teissier's result for hypersurfaces of $\C^3$ which establishes the 
 Whitney equisingularity of $X$ and its tangent cone under this conditions.  
\end{abstract}

\section{Introduction}
 
    The goal of this paper is to take a step in the study of the geometry of the
  specialization space $\varphi: (\X,0) \to (\C,0)$ of a germ of reduced and  $d$ dimensional 
  singularity $(X,0)$ to its tangent cone $C_{X,0}$ from the point of view of Whitney \linebreak
  equisingularity. The map $\varphi$ describes a flat family of analytic germs with a section
  $ \X \stackrel{\curvearrowleft}{\to} \C :\sigma$, such that for each $t \in \C^*$ the germ $(\varphi^{-1}(t),\sigma(t))$
  is isomorphic to $(X,0)$ and the special fiber is isomorphic to the tangent cone. This construction is 
  essentially due to Gerstenhaber \cite{Ger} in a more algebraic setting.\\

   One would like to establish conditions on the strata of the canonical Whitney stratification
 of a reduced complex analytic germ which ensure the Whitney equisingularity of the germ and
 its tangent cone. In this paper we achieve the "codimension zero" part of this program.\\
  
    The space $(\X,0)\to (\C,0)$ has been used to study Whitney conditions in \cite{Na}, and to study
  the structure of the set of limits of tangent spaces in \cite{L-T2} and \cite{L-T1}. In \cite{L-T2}, 
  the authors prove the existence of a finite family $\{V_\alpha\}$ of subcones of the reduced tangent
  cone $|C_{X,0}|$  that determines the set of limits of tangent spaces to $X$ at $0$.\\
  
  	To be more specific, we fix an embedding $(X,0)\subset (\C^{n+1},0)$ and build the normal/conormal 
   diagram, 
  	\[\xymatrix{E_0C(X)\ar[r]^{\hat{e}_0}\ar[dd]^{\kappa'}\ar[ddr]^\xi & C(X)\ar[dd]^\kappa
             \\
             & &  \\
             E_0X\ar[r]_{e_0}  &  X }\]
    where $E_0X \subset X \times \P^n$ is the blowup of $X$ at the origin, $C(X) \subset X \times
   \check{\P}^n$ is the conormal space of $X$ whose fiber determines the set of limits of tangent 
   spaces (see section \ref{conormal}), and $E_0C(X) \subset X \times \P^n 
   \times \check{\P}^n$ is the blowup in $C(X)$ of the subspace $\kappa^{-1}(0)$; consider the irreducible
   decomposition of the reduced fiber $|\xi^{-1}(0)|= \bigcup D_\alpha$. The authors prove that
   the fiber $\xi^{-1}(0)$ is contained in the incidence variety $I \subset \P^n \times \check{\P}^n$
   and that each $D_\alpha$ establishes a projective duality of its images $V_\alpha \subset \P 
   C_{X,0} \subset \P^n$ and $W_\alpha \subset \kappa^{-1}(0) \subset \check{\P}^{n}$.\\

   In particular, the $V_\alpha$'s that are not irreducible components of the tangent cone are called 
  exceptional cones and they appear in $\X$ as an obstruction to the $a_f$ stratification of the 
  morphism $\X \to \C$. They also prove that if the germ $(X,0)$ is a cone itself, then 
  it doesn't have exceptional cones. So a natural question arises, if a germ of analytic singularity
  $(X,0)$ doesn't have exceptional tangents, how close is it to being a cone?\\
  
    A partial answer to this question was given in \cite{L-T1} in terms of Whitney equisingularity.
   The authors prove that for a surface $(S,0) \subset (\C^3,0)$ with reduced tangent cone $C_{S,0}$, 
   the absence of exceptional cones is a necessary and sufficient condition for it to be Whitney 
   equisingular to its tangent cone. \\
  
     The specialization space $(\X,0) \to (\C,0)$ has a canonical section which picks the  
  origin in each fiber (see section 2). Let $Y \subset \X$ be given by this section and let $\X^0$ 
  be the non singular part $\X$. The main objective of this paper is to prove that if the germ $(X,0)$
  doesn't have exceptional cones and the tangent cone is reduced, then the couple $(\X^0, Y)$ satisfies 
  Whitney's conditions a) and b) at the origin. 
    
\section{Specialization to the tangent cone.}

   	Let $(X,0)$ be a reduced germ of analytic singularity of pure dimension $d$, with tangent cone
$C_{X,0}$. Recall that the projectivized tangent cone can be defined as the exceptional divisor of the 
blowup of $X$ in $0$, and it is equivalent to considering the analytic ``proj'' of the graded algebra
     \[gr_{\m}O_{X,0}\colon = \bigoplus_{i\geq0}\m^i/ \m^{i+1}\]
where $\m$ is the maximal ideal of the analytic algebra $O_{X,0}$ associated to the germ. Moreover,
if we consider an embedding $(X,0) \subset (\C^{n+1},0)$, the analytic algebra $O_{X,0}$ is isomorphic to
$\C\{z_0,\ldots,z_n\}/ I$, where $I$ is an ideal, $gr_{\m}O_{X,0}$ is isomorphic to $\C[z_0,\ldots,z_n]/\textrm{In}_{\M}I$ where
$\M$ is the maximal ideal of $\C\{z_0,\ldots,z_n\}$, and the ideal In$_{\M}I$ is generated by all the initial forms
with respect to the $\M$-adic filtration of elements of $I$.\\

  Let us suppose that the generators $\left<f_1,\ldots,f_p\right>$ for $I$, were chosen in such a way that their 
 initial forms generate the ideal In$_{\M}I$ defining the tangent cone. Note that the $f_i$'s are convergent
 power series in $\C^{n+1}$, so if $m_i$ denotes the degree of the initial form of $f_i$, by defining 
   \begin{equation}
 \label{eq1} F_i(z_0,\ldots,z_n,t):= t^{-m_i}f_i(tz_0,\ldots,tz_n)
 \end{equation} 
 we obtain convergent power series, defining holomorphic functions on a suitable open subset $U$ of 
 $\C^{n+1} \times \C$. Moreover, we can define the analytic algebra
 \[O_{\X,0}= \C \{z_0, \ldots, z_n,t\}/\left< F_1, \ldots, F_p\right>\]
 with a canonical morphism $\C\{t\} \longrightarrow O_{\X,0}$ coming from the inclusion
 $\C\{t\} \hookrightarrow \C\{z_0, \ldots, z_n,t\}$. Corresponding to this morphism of analytic algebras,
 we have the map germ $\varphi: (\X,0) \to (\C,0)$ induced by the projection of $\C^{n+1} \times \C$ to 
 the second factor.

\begin{definition}
    The germ of analytic space over $\C$, 
   \[\varphi: (\X,0) \to (\C,0) \]
    is called the specialization of $(X,0)$ to its tangent cone $(C_{X,0},0)$.
\end{definition}
 
   There is another way of building this space that will allow us to derive some interesting properties. Let
$E_{(0,0)}\C^{n+2}$ be the blowing up of the origin of $\C^{n+2}$, where we now have the 
coordinate system $(z_0, \ldots, z_n , t)$. Let $W \subset E_{(0,0)}\C^{n+2}$ be the chart where the 
invertible ideal defining the exceptional divisor is generated by $t$, that is, in this chart the blowing up 
map is given by $(z_0, \ldots,z_n,t) \mapsto (tz_0, \ldots, tz_n,t)$.
\[ \xymatrix { W \ar @{^{(}->}[r] \ar[dr] & E_{(0,0)}\C^{n+2} \ar[d]^{E_0}\\
                 & \C^{n+2} }\]
 
\begin{lema}\label{Specialisationasblowup}
 Let  $X \times \C \subset \C^{n+2}$ be a small enough representative of the germ 
$(X \times \C,0)$. If $(X \times \C)'$ denotes the strict transform of $(X \times \C)$ in the blowing up 
$E_{(0,0)}\C^{n+2}$, then the space $(X \times \C)' \bigcap W$ together with the map induced by the 
restriction of the map $E_{(0,0)}\C^{n+2} \to \C^{n+1} \times \C \to \C$ is isomorphic to the 
specialization space $\varphi: \X \to \C$.
\end{lema}
\begin{proof}
   We know that the strict transform $(X \times \C)'$ is isomorphic to the blowing up of $X \times \C$
at the origin, and we are seeing it as a reduced analytic subvariety of $\C^{n+2} \times \P^{n+1}$. 
This means that the exceptional divisor $(X \times \C)' \cap (\{0\}\times \P^{n+1})$ is equal
to $\P(C_{X,0} \times \C)$, and so the ideal defining it is generated by the ideal defining
the tangent cone $C_{X,0}$ in $\C^{n+1}$, that is, the ideal of initial forms In$_{\M}I$.
By hypothesis, $W \subset E_{(0,0)}\C^{n+2} \subset \C^{n+2} \times \P^{n+1}$ is set theoretically
described by 
\[W=\left\{(tz_0,\ldots,tz_n,t), [z_0:\cdots:z_n:1] \, | \, (z_0,\ldots,z_n,t) \in \C^{n+2}\right\}\]
so in local coordinates the map $E_0$ restricted to $W$ is given by $(z_0, \ldots,z_n,t) \mapsto (tz_0,
 \ldots,tz_n, t)$. Finally, since the ideal defining $X \times \C$ is generated in $\C\{z_0,\ldots,z_n,t\}$ 
by the ideal $I=\left<f_1,\ldots,f_p\right>$ of $\C\{z_0,\ldots,z_n\}$ defining $X$ in $\C^{n+1}$, and since 
we have chosen the $f_i$'s in such a way that their initial forms generate the ideal In$_{\M}I$, then
the ideal defining the strict transform $(X \times \C)'$ in $W$ is given by
\[\mathfrak{J}O_W=\left< t^{-m_1}f_1(tz_0, \ldots, tz_n), \ldots, t^{-m_p}f_p(tz_0,\ldots,tz_n)\right>O_W\]
that is, we find the same functions $F_1, \ldots, F_p$ which we used to define $\varphi:\X \to \C$.
\end{proof}

 \begin{proposition}  \label{PropertiesSpecialisation}
 Let $\varphi:\X \to \C$ be a small enough representative of the germ, then:
    \begin{enumerate}
 \item The morphism $\varphi$ is induced by the restriction of the projection \linebreak
       $\C^{n+1} \times \C \to \C$ to the closed subspace defined by $(F_1,\ldots,F_p)$, and it is faithfully flat.
 \item The special fiber $\X(0)\colon = \varphi^{-1}(0)$ is isomorphic to the tangent cone $C_{X,0}$.
 \item The analytic space $\X \setminus \varphi^{-1}(0)$ is isomorphic to $X \times \C^*$ as an analytic 
       space over $\C^*$. In particular, for every $t \in \C^*$, the germ $(\varphi^{-1}(t),\{0\} \times t)$ is 
       isomorphic to $(X,0)$. 
 \item The germ $(\X,0)$ is reduced and of pure dimension $d+1$.
\end{enumerate}
  that is, we have produced a 1-parameter flat family of germs of analytic spaces
specializing $(X,0)$ to $(C_{X,0},0)$.
 \end{proposition}
\begin{proof}
   First of all, note that the inclusion $\C\{t\} \hookrightarrow \C\{z_0, \ldots, z_n,t\}$ can be seen
as the stalk map at the origin of the holomorphic map defined by the linear projection onto the last coordinate
$\C^{n+1} \times \C \to \C$. This implies that $\varphi$ is just the restriction to $\X$ of this projection.

   Now, to prove the (faithful) flatness of $\varphi$ we must prove that $O_{\X,0}$ is faithfully flat 
as a $\C\{t\}$ module, but by \cite[Prop. B.3.3, p. 404]{G-L-S} flat implies faithfully flat for local rings, and
by \cite[Corollary 7.3.5, p. 390 ]{Gre-Pfis} $O_{\X,0}$ is flat if and only if it is torsion free. In other words
all we have to prove is that $t$ is not a zero divisor in $O_{\X,0}$.

   But by lemma \ref{Specialisationasblowup}, $\X$ is isomorphic to an open subset of the blowing up of 
$X \times \C$ along the subspace $\{0\} \times \C$, where the ideal of the exceptional divisor is invertible, 
generated by $t$. Thus, by definition of blowing up, $t$ is not a zero divisor, $\X$ is of pure dimension $d+1$ 
(the dimension of $X \times \C$), and since the blowing up of a reduced space remains reduced then $\X$ is
reduced.\\ 

%The germ $(\X,0)$ is reduced (\textbf{WHY???}), so it has no embedded componentns and all the associated 
%prime ideals $\{P_1, \ldots, P_r\}$ of the zero ideal $\left<0\right>O_{\X,0}$ are minimal. If $t$ is a zero divisor
%then there exists $j \in \{1, \ldots,r\}$ such that $t \in P_j$ and so the irreducible component $\X_j$ 
%of $\X$ corresponding to this prime ideal is contained in the hyperplane $\{t=0\}$ of $\C^{n+2}$. This means
%by \cite[Prop. 49.5, p. 193]{Ka} that the closure of $\X^\circ \setminus \X_j$ does not contain $\X_j$,
%and the same goes for $\X^\circ \setminus \{t=0\}$. But this is a contradiction with the way $\X$ was 
%constructed as a strict transform in lemma \ref{Specialisationasblowup}.  

%   Moreover, by definition, $\X$ is the zero locus of the analytic functions:
%\[F_j(\underline{z},t)= t^{-m_j}f_j(t\underline{z})=f_{m_j}+tf_{m_j+1}+ t^2f_{m_j+2} + \cdots \]
%where $f_{m_j}$ is the initial form of $f_j$ and has degree $m_j$. So, the equations defining the special 
%fiber $\X(0)$ are $\{t=0, f_{m_1}=0, \ldots, f_{m_p}=0\}$ and by hypothesis these initial forms generate
%the ideal of initial forms In$_{\M}I$ defining the tangent cone $C_{X,0}$.

   The biholomorphism of the map induced by the isomorphism \[\phi: \C^{n+1} \times \C^* \to \C^{n+1} \times \C^*
\mathrm{defined \; by} (\underline{z},t) \mapsto (t\underline{z},t)\] is also a direct consequence of lemma 
\ref{Specialisationasblowup}. It maps $\X \setminus \X(0)$ onto $\X \times \C^*$, and for each $t\neq 0$ the 
fiber $\X(t)$ is mapped bilohomorphically onto $X \times \{t\}$.
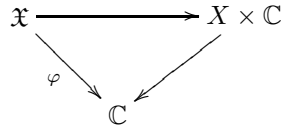
\[\xymatrix { \X \ar[dr]_{\varphi} \ar[rr]  &  & X \times \C \ar[dl] \\
               & \C & }\] 
 
  Finally, the fact that the special fiber $\X(0)$ is isomorphic to the tangent cone can be read directly 
from the analytic functions $F_1, \ldots, F_p$ defining $\X$, since when setting $t=0$ we have the initial 
forms $F_i(z,0)=f_{m_i}$ which by hypothesis generate the ideal defining the tangent cone in 
$\C^{n+1}$.

\end{proof}

   A more detailed description of this space, relating it to a generalized Rees algebra and  
 interpreting the space thus obtained as the open set of the blowup (\ref{Specialisationasblowup}) 
 of $X \times \C$ at the origin can be found in \cite[p. 428-430]{L-T1} for surfaces, and 
 \cite[p. 556-557]{L-T2}, or \cite[p. 200-202]{Na} in the general case.

\begin{remark} Note that:
   \begin{enumerate}
    \item  The map $\phi: \X \to X \times \C$ from proposition \ref{PropertiesSpecialisation}
    is defined everywhere and maps the entire fiber $\X(0)$ to the origin in $X \times \C$.
    \item If we denote by $\X(t)^0$ the non-singular part of the fiber, the open dense subset 
     $\bigcup_{t} \X(t)^0 \subset \X$ is called the \textbf{relative smooth locus of $\X$ with respect
     to $\varphi$}.
   \end{enumerate}
\end{remark}     

\begin{lema} \label{Especializaciones}
     Let $X=\cup_{j=1}^r X_j$ be the irreducible decomposition of $X$. Then, 
  the specialization $\X_j$ of $X_j$ is an analytic subspace of $\X$, and the following diagram
  commutes. 
   
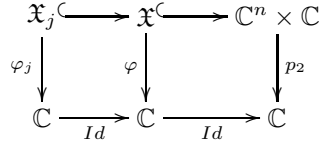
\[ \xymatrix { \X_j \ar@{^{(}->}[r]\ar[d]_{\varphi_j} &  \X \ar@{^{(}->}[r]\ar[d]_{\varphi} & \C^n \times \C \ar[d]^{p_2}  \\
                   \C \ar[r]_{Id} & \C \ar[r]_{Id} & \C } \]
  In particular, $\X=\cup_{j=1}^r \X_j$ is the irreducible decomposition of $\X$.
  \end{lema}
  \begin{proof}$\;$\\
  
        Note that $X_j$ is a proper analytic subspace of $X$ for all $j\geq 1$, so we have a strict inclusion
    of their corresponding ideals in $O_{n+1}:=\C\{z_0,\ldots,z_n\}$, namely $I \subset J$, from which we immediately
    obtain that $In_\M I \subset In_\M J$ or equivalently $C_{X_j,0} \subset C_{X,0}$. \\
    
        Now let us take as before, generators for $I$, say $I=\left<f_1,\ldots,f_p\right>$, in such a way that 
    their initial forms generate the ideal defining the tangent cone $In_{\M}I$, and doing the same for $J$,
    we get $J=\left<g_1,\ldots,g_s\right>$ and $In_\M J= \left<in_\M g_1,\ldots, in_\M g_s\right>$. But the previous
    inclusions tell us that we can choose as generators for $J=\left<f_1,\ldots,f_p,g_1,\ldots,g_s\right>$, and still
    get that their initial forms generate the ideal $In_\M J =\left<in_\M f_i,in_\M g_j\right>$.\\
    
        So finally, to build the specialization spaces $\X$ and $\X_j$ as we did before, we define the convergent series
    in $O_{n+2}$, $F_i(z,t)=t^{-m_i}f_i(tz_0,\ldots,tz_n)$ and $G_j(z,t)=t^{-m_j}g_j(tz_0,\ldots,tz_n)$, that give
    us the embedding $\X:= V(F_1,\ldots,F_p) \subset \C^{n+1} \times \C$ and the embedding $\X_j := V(F_1,\ldots,F_p,G_1,
    \ldots,G_s) \subset \C^{n+1} \times \C$. Moreover, since $\left < F_1, \ldots,F_p\right> \subset 
    \left < F_i, G_j\right>$ then we have a closed embedding $\X_j \subset \X$ compa- tible with the projection to the 
    $t$ axis.\\
    
    	   And even more, since with respect to this embedding of $\X$ in $\C^{n+1} \times \C$, the isomorphism 
    $\phi$ is of the form:
    \begin{align*}
              \phi: \X \setminus \varphi^{-1}(0) & \longrightarrow X \times \C^*\\
                     (z_0,\ldots,z_n,t)    & \longmapsto (tz_0,\ldots,tz_n,t)
     \end{align*}
    We also have compatibility with the isomorphism, that is $\phi_j = \phi|_{\X_j}$.
    \[ \xymatrix { \X \setminus \varphi^{-1}(0)\ar [r]^\phi & X \times \C^*  \\
                   \X_j\setminus \varphi_j^{-1}(0) \ar@{^{(}->}[u]\ar[r]_{\phi_j} & X_j \times \C^* \ar@{^{(}->}[u]} \]
  \end{proof}

\begin{remark}\label{SpecialisationNormalCone}
 \begin{enumerate}
  \item For an analytic subspace $Y \subset X$ we can mimic the cons- truction of lemma 
  \ref{Specialisationasblowup} to build the specialization 
  space $\varphi: \X \to \C$ where we still have that the fiber $(\X(t),(0,t))_{t\neq 0}$ is isomorphic
  to the germ $(X,0)$, but this time the special fiber $(\X(0),(0,0))$ is isomorphic to the
  normal cone $(C_{X,Y},0)$. The map $\varphi$ is again faithfully flat.
   \item If  $Y \subset X$ is a linear subspace defined by the ideal \linebreak
         $J=\left<z_0, \ldots z_{n-s} \right>\C\{z_0,\ldots,z_{n-s},y_1,\ldots,y_s\}$ then we can choose
    analytic functions $f_1, \ldots, f_p$ such that they generate the ideal $I$ defining $X$ in $\C^{n+1}$, and
    their initial forms $f_{m_i}=in_J f_i$ generate the ideal of initial forms $in_J I$. In this case the  
    ideal generated by the analytic functions $F_i(z,y,t)= t^{-m_i}f_i(tz_0,\ldots,tz_{n-t},y_1,\ldots,y_s)$
    will be the ideal defining the space $\X$  in $\C^{n+1} \times \C$, where $m_i$ is equal to
    $\nu_Y f_i$. 
         
 \end{enumerate}
\end{remark}

\section{The Relative Nash Modification of $\mathbf{\X}$}\label{LimitsofTanSpacesSection}

   Let us take a representative $\varphi: \X \to \C$ of the germ $\varphi:(\X,0) \to (\C,0)$, and consider
the map:
\begin{align*} \gamma_\varphi: \X^\circ_\varphi & \longrightarrow Gr(d,n+1)\\
                              (z,t)  & \longrightarrow T_{(z,t)}X^\circ_\varphi(t)
\end{align*} 
where $\X^\circ_\varphi$ denotes the relative smooth locus of $\X$ with respect to $\varphi$, 
$Gr(d,n+1)$ corresponds to the grassmannian of directions of $d-$planes of the hyperplane $\{t=0\} \subset 
\C^{n+1} \times \C$, and  $T_{(z,t)}\X^\circ_\varphi(t)$ denotes the tangent space to the fiber
$\X^\circ_\varphi(t)$ at the point $(z,t)$. The closure $\Na_\varphi \X$ of the graph of $\gamma_\varphi$ in 
$\X \times Gr(d,n+1)$ is an analytic space of dimension $d+1$, which is known as \textbf{the relative Nash
modification of $\varphi:\X \to \C$.}\\

    A. Nobile proved in \cite[Thm 1, p. 299]{No} that the Nash modification is a blowing up.
The main ingredient of his proof is the Plucker embedding of the grassmannian $G(d,n+1)$ in the 
projective space $\P^{N}$, where $N=\begin{pmatrix} n+1 \\ d \end{pmatrix}$. Minor modifications of the proof 
immediately gives us an analogous result for the relative case. We will only state it in the case of 
$\varphi: \X \to \C$.

\begin{lema}\label{NashNobileCoro}
  The relative Nash modification $\nu_\varphi:\Na_\varphi \X \to \X$ is locally a blowing-up with center a 
suitable ideal $J_\varphi \subset O_{\X}$. Moreover, 
if $(\X,0)$ is a complete intersection of dimension $n+2-p$ then we may take the ideal $J_\varphi \subset
O_\X$ to be the relative Jacobian ideal, formed by the $ p \times p$ minors of the relative Jacobian matrix 
$[D_\varphi F]=\left[\frac{\partial F_i}{\partial z_j}\right]^{i=1\ldots p}_{j=0 \ldots n}$. (We are 
omitting the partial derivatives with respect to the parameters, which in this case correspond
to the t-coordinate).
\end{lema}
\begin{proof}
    Given integers $n+1 \geq r > 0, \; p \geq n+1-r$ and a $p \times (n+1)$ matrix $A$, let 
$S$(resp. $S'$) denote the set of increasing sequences of $n+1-r$-positive integers less than $p+1$ 
(resp. $n+2$); if $\alpha=(\alpha_1, \ldots, \alpha_{n+1-r}) \in S$, $\beta=(\beta_1, \ldots, \beta_{n+1-r})\in S'$,
then $M_{\alpha \beta}$ will denote the minor of $A$ obtained by considering the rows determined by $\alpha$
and the columns determined by $\beta$.

   Following the proof of Nobile, let $\X=\bigcup_{j=1}^k \X_j$ be the irreducible decomposition
of a small enough representative of $(\X,0)$. Let $[D_\varphi F]=
\left[\frac{\partial F_i}{\partial z_j}\right]^{i=1\ldots p}_{j=0 \ldots n}$ be the relative Jacobian
matrix of the map $\varphi: \X \to \C$. By construction, there is an open dense set $U \subset \X$,
such that for every point $(z_0,t_0)$ in $U$ the matrix $[D_\varphi F(z_0,t_0)]$ has rank $n+1-d$. Since
$\X$ is reduced, each irreducible component $\X_i$ is reduced and so for each $i= 1, \ldots, k$ there exists 
a pair $(\alpha^i, \beta^i) \in S \times S'$ such that the $(n+1-d) \times (n+1-d)$ minor $M_{\alpha^i \beta^i}$ 
of $[D_\varphi F]$ does not vanish identically on $\X_i$. For each $i= 1, \ldots, k$, fix $H_i \in O_{\X,0}$
such that $H_i=0$ on $\bigcup_{j\neq i} \X_j$, and $H_i \neq 0$ on $\X_i$. For each $\beta \in S'$ define 
the function $G_\beta= \sum_{i=1}^k H_iM_{\alpha^i \beta} \in O_{\X,0}$, and consider the ideal 
$J_\varphi \subset O_{\X,0}$ generated by the $G_\beta$'s.  

   Note that the analytic subset $V(J_\varphi)$ of $\X$ defined by the ideal $J_\varphi$ contains the 
relative singular locus of $\varphi:\X \to \C$. Moreover, the open set $W:=\X \setminus V(J_\varphi)$ is dense
in $\X$. Finally if we build a representative of this blowup using the functions $G_\beta$, we will have it 
as an analytic subspace of $\X \times \P^N$, with $N=\begin{pmatrix} n+1 \\ n+1-d \end{pmatrix}-1= \begin{pmatrix} 
n+1 \\ d \end{pmatrix}-1$, and for a point $(z,t) \in \X_i \cap W$ we have that:
\[[G_\beta(z,t)]=[\sum_{j=1}^k H_j(z,t)M_{\alpha^j \beta}(z,t)]=[H_i(z,t) M_{\alpha^i \beta}(z,t)]=
   [M_{\alpha^i \beta}(z,t)] \in \P^N \]
which corresponds to the coordinates of the tangent space $T_{(z,t)}X^\circ_\varphi(t)$ for the Plucker
embedding of the grassmannian $G(d,n+1)$, in the projective space $\P^{N}$.
\end{proof}

  This lemma allows us to establish the following relation between the Nash modification of $X$ and
 the relative Nash modification of $\X$.

\begin{proposition}\label{NashrelativevsNash}
   There exists a natural surjective morphism $\Gamma: \Na_\varphi \X \to \Na X $, making the following 
diagram commute:
\[\xymatrix{ \Na_\varphi \X \ar[r]^\Gamma \ar[d]_{\nu_\varphi} &  \Na X \ar[d]^{\nu}\\
                   \X  \ar[r]_\phi & X }\]   
\end{proposition}
\begin{proof}

   Algebraically, this results from the universal property of the blowup \linebreak $\nu: \Na X \to X$. We
start with the diagram:
\[\xymatrix{ \Na_\varphi \X \ar[d]_{\nu_\varphi} &  \Na X \ar[d]^{\nu}\\
                   \X  \ar[r]_\phi & X }\]   
where the map $\phi$ is defined by $(z_0,\ldots,z_n,t) \to (tz_0, \ldots, tz_n)$, and so it induces
a morphism of analytic algebras $\phi^*: O_{X,0} \to O_{\X,0}$ defined by $z_i \to tz_i$.\\

 Recall that the ideal of the germ $(\X,0)$ is generated by the series $F_i(z,t)= t^{-m_i}f_i(tz) \in 
\C\{z_0, \ldots, z_n,t\}$, $i=1, \ldots, p$, where the series $f_j \in \C\{z_0, \ldots, z_n\}$ are such that 
they generate the ideal of $(X,0)$ in $(\C^{n+1},0)$ and their initial forms generate the ideal of 
$(C_{X,0},0)$.\\  

  By \cite[Thm 1, p. 299]{No} there exists an ideal $JO_{X,0}$ whose blowup is isomorphic to the Nash modification 
of $X$. We have to prove that the ideal $\phi^*(J)O_{\X,0}$ is locally invertible when pulled back 
to $\Na_\varphi \X$.\\
  
   Let $X= \bigcup_{j=1}^k X_j$ be the irreducible decomposition of a small enough representative of $(X,0)$.
Then the irreducible decomposition of a small enough representative of the germ $(\X,0)$ is of the form
$\bigcup_{j=1}^k \X_j$, where for each $j$ the space $\X_j$ is isomorphic to the specialization space of the
$X_j$ component to its tangent cone $C_{X_j,0}$.  Now, by \cite[Thm 1, p. 299]{No} the ideal $J\subset O_{X,0}$ can be constructed
in the following way (see the proof of \ref{NashNobileCoro} for more details and notation): 
For each $i=1, \ldots, k$ there exists a pair $(\alpha^i, \beta^i) \in S \times S'$ such that 
the $(n+1-d) \times (n+1-d)$ minor $\mu_{\alpha^i \beta^i}$ 
of the jacobian matrix $[Df]$ does not vanish identically on $X_i$. Then for each $i= 1, \ldots, k$, 
choose a function $h_i \in O_{X,0}$ such that $h_i=0$ on $\bigcup_{j\neq i} X_j$, and $h_i \neq 0$ on $X_i$.
By taking powers of the $h_i$'s if necessary we can assume they are all of the same order $\gamma$. 
Finally for each $\beta \in S'$ define the function $g_\beta= \sum_{i=1}^k h_i\mu_{\alpha^i \beta} \in O_{X,0}$, and 
define $J$ as the ideal generated by the $g_\beta$'s.\\

   Consider an $(n+1-d) \times (n+1-d)$ minor $\mu_{\alpha \beta}$ of the jacobian matrix $[Df]$
\[\mu_{\alpha,\beta}=\begin{vmatrix} \frac{\partial f_{\alpha_1}}{\partial z_{\beta_1}}(z) & \cdots & 
                       \frac{\partial f_{\alpha_1}}{\partial z_{\beta_{n+1-d}}}(z) \\
                           \vdots            &  \vdots &    \vdots \\
                      \frac{\partial f_{\alpha_{n+1-d}}}{\partial z_{\beta_1}}(z) & \cdots & 
                       \frac{\partial f_{\alpha_{n+1-d}}}{\partial z_{\beta_{n+1-d}}}(z)
      \end{vmatrix} \]
  Then, from the equalities $\phi^*(\frac{\partial f_i}{\partial z_j}(z))= 
\frac{\partial f_i}{\partial z_j}(tz)$, and $\frac{\partial f_i}{\partial z_j}(tz)= t^{m_i-1} \frac{\partial F_i}
{\partial z_j}(z,t)$, we have that the minor $\mu_{\alpha \beta}$ is mapped under $\phi^*$ to:
\begin{align*}
 \phi^*(\mu_{\alpha \beta}) & =\begin{vmatrix} 
  t^{m_{\alpha_1}-1}\frac{\partial F_{\alpha_1}}{\partial z_{\beta_1}}(z,t) & \cdots & t^{m_{\alpha_1}-1}
                  \frac{\partial F_{\alpha_1}}{\partial z_{\beta_{n+1-d}}}(z,t) \\
                \vdots            &  \vdots &    \vdots \\
   t^{m_{\alpha_{n+1-d}}-1}\frac{\partial F_{\alpha_{n+1-d}}}{\partial z_{\beta_1}}(z,t) & \cdots & 
   t^{m_{\alpha_{n+1-d}}-1}\frac{\partial F_{\alpha_{n+1-d}}}{\partial z_{\beta{n+1-d}}}(z,t)
      \end{vmatrix}\\
      & =
 t^{(\sum_1^{n+1-d} m_{\alpha_i})-(n+1-d)}\begin{vmatrix} 
    \frac{\partial F_{\alpha_1}}{\partial z_{\beta_1}}(z,t) & \cdots & 
                  \frac{\partial F_{\alpha_1}}{\partial z_{\beta_{n+1-d}}}(z,t) \\
                \vdots            &  \vdots &    \vdots \\
    \frac{\partial F_{\alpha_{n+1-d}}}{\partial z_{\beta_1}}(z,t) & \cdots & 
    \frac{\partial F_{\alpha_{n+1-d}}}{\partial z_{\beta{n+1-d}}}(z,t)
      \end{vmatrix}\\
      & = t^{(\sum_1^{n+1-d} m_{\alpha_i})-(n+1-d)} M_{\alpha \beta}
\end{align*}
  where $M_{\alpha \beta}$ is the $(n+1-d) \times (n+1-d)$ minor of the relative jacobian matrix
 $[D_\varphi F]$.\\

   If we define $H_i \in O_{\X,0}$ by $H_i(z,t)=t^{-\gamma}h_i(tz)$, then each $H_i$ satisfies
that $H_i=0$ on $\bigcup_{j\neq i} \X_j$, and $H_i \neq 0$ on $\X_i$ and so for each $\beta \in S'$
we have that
\[\phi^*(g_\beta)= \sum_{i=1}^k \phi^*(h_i)\phi^*(\mu_{\alpha^i \beta}) =
           t^{\left(\gamma +(\sum_1^{n+1-d} m_{\alpha_i})-(n+1-d)\right) }\sum_{i=1}^k H_iM_{\alpha^i \beta}
  = t^r G_\beta \]
and so 
\[\phi^*(J)O_{\X,0}=\left< t^r\right>J_\varphi O_{\X,0}\]
where by the proof of \ref{NashNobileCoro} $J_\varphi O_{\X,0}$ is an ideal whose blowup is isomorphic
to the relative Nash modification $\Na_\varphi \X$. But by definition of the blowup, the ideal 
$J_\varphi O_{\X,0}$ is locally invertible when pulled back to $\Na_\varphi \X$. It follows that after 
multiplication  by the invertible ideal $\left<t^r\right>$ in $O_{\X,0}$, it will remain locally invertible when pulled back to 
$\Na_\varphi \X$.\\

  Finally, note that for the diagram to be commutative, the morphism $\Gamma$ must map the point
$(z,t,T_{(z,t)}\X^\circ_\varphi(t)) \in \Na_\varphi \X$ to the point $(tz, T_{(tz)}X^\circ) \in \Na X$.
That is the tangent space $T_{(z,t)}\X^\circ_\varphi(t)$ to the 
fiber $\X(t)$ is canonically identified with the tangent space $T_{(tz)}X^\circ$ to $X$ at the corresponding
points. As it should be since we know that the restriction of the map $\phi$ to any fiber 
$(\X(t),(0,t))$ for $t\neq 0$ is an isomorphism with $(X,0)$.  
\end{proof}

 \section{The conormal space and relative conormal space of $\X$.}\label{conormal}
  
	   Let $\X\subset \C^{n+2}$ be a representative of the germ $(\X,0)$. Recall that the
 projectivized conormal space of $\X$ in $\C^{n+2}$ is an analytic space $C(\X) \subset \X \times \check{\P}^{n+1}$,
	together with a proper analytic map $\kappa_\X: C(\X)\to \X$, where the fiber over a smooth point 
	$p \in \X$ is the set of tangent hyperplanes, that is the hyperplanes $H$ containing the 
  direction of the tangent space $T_p\X$. The space $C(\X)$, depends on the embedding, however
  the fiber $\kappa_\X^{-1}(p)$ allows us to recover the fiber of the Nash modification, which 
  is independent of the embedding. Up to now  we have:
\[\Na \X \subset \X \times G(d+1,n+2) \subset \X \times \P^N\]
But we know that the grassmannian $G(d+1,n+2)$ is isomorphic to the grassmannian $G(n+1-d,n+2)$ and
the isomorphism is given by sending a $d+1-$plane $T$ to the $n+1-d-$plane $L$ of linear functionals
in $\check{\C}^{n+2}$ that vanish on $T$. With this isomorphism, we have:
\[\Na \X \subset \X \times G(n+1-d,n+2) \subset \X \times \P^N \]  
  Let $\Xi \subset G(n+1-d, n+2) \times \check{\P}^{n+1}$ denote the tautological bundle, that is
$\Xi = \{(L,[a])\; | \; L \in G(n+1-d,n+2), \; [a] \in \P L \subset \check{\P}^{n+1} \}$, and consider
the intersection 
\[ \xymatrix{ 
  E:=\{\X \times \Xi\} \bigcap \{\Na \X \times \check{\P}^{n+1} \} \ar @{^{(}->}[r] \ar[dr]^{p_2} \ar[d]_{p_1}&  
  \X \times G(n+1-d,n+2) \times \check{\P}^{n+1} \ar[d]\\
   \Na \X & \X \times \check{\P}^n}\]
with the vertical morphism $p_2$ being the morphism induced by the projection onto $ \X \times \check{\P}^{n+1}$.
We then have the following result.

\begin{proposition}\label{ConormalvsNash}
   Let $p_2:E \to \X \times \check{\P}^{n+1}$ be as before. The set theoretical image $p_2(E)$ of the morphism 
$p_2$ coincides with the conormal space of $\X$ in $\C^{n+2}$ 
\[ C(\X) \subset \X \times \check{\P}^{n+1}\] 
Moreover, the morphism $p_1: E \to \Na \X$ is a locally trivial fiber bundle over $\nu^{-1}(\X^0) \subset
 \Na \X$ with fiber $\P^{n-d}$.
%Moreover, the morphism $p_2:E \to C(X)$ is finite over $C(X^0)$%
\end{proposition}
\begin{proof}
  By definition, the conormal space of $\X$ in $\C^{n+2}$ is an analytic space $C(\X) \subset X \times 
  \check{\P}^{n+1}$, together with a proper analytic map $\kappa: C(\X)\to \X$, where the fiber over a smooth point 
  $x \in \X^0$ is the set of tangent hyperplanes, that is the hyperplanes $H$ containing the direction of the 
  tangent space $T_x\X$. That is, if we define $E^0=\{(x,T,[a]) \in E \, | \, x \in \X^0\}$, then by 
  construction $E^0= p_1^{-1}(\nu^{-1}(\X^0))$, and $p_2(E^0)= C(\X^0)$. Since the morphism $p_2$ is proper, in 
  particular it is closed which finishes the proof.
\end{proof} 

  In the same way we can construct
  the relative conormal space $C_\varphi(\X)$ as a subvariety of $\X \times \check{\P}^n$ where 
  $\check{\P}^n$ stands for the dual projective space of directions of hyperplanes of the hyperplane
  $\{t=0\} \subset \C^{n+1} \times \C$ .\\
  
 \begin{proposition}\label{IsomorfismosConormales}
    Let $Y \subset X$ be a smooth analytic subvariety of dimension $0 \leq s < d$, let 
$\varphi:\X \to \C$ denote the specialization space of $X$ to its normal cone along $Y$, and
let $\phi: \X \to X \times \C$ denote the canonical map obtained from the construction in 
lemma \ref{Specialisationasblowup}. Then there exist isomorphisms $\psi: C(\X \setminus \X(0)) \to C(X) \times \C^*$;
$P:C(\X \setminus \X(0)) \to C_\varphi(\X \setminus \X(0))$; and $\psi_\varphi: C_\varphi(\X \setminus \X(0)) \to C(X) \times \C^*$ 
making the following diagram commutative:

  \[ \xymatrix  {  C(\X \setminus \X(0)) \ar[d]^P \ar[r]^{\psi}  & C(X) \times \C^* 
              \ar[d]^{Id} \ar[r]^{\widetilde{pr_1}} & C(X) \ar[d]^{Id}\\
              C_\varphi(\X \setminus \X(0)) \ar[d]^{\kappa_\varphi} \ar[r]^{\psi_\varphi}  & C(X) \times \C^* 
              \ar[d]^{\kappa_X \times Id} \ar[r]^{\widetilde{pr_1}} & C(X) \ar[d]^{\kappa_X}\\
                  \X \setminus \X(0) \ar[r]_{\phi} \ar[dr]_\varphi &  X \times \C^* \ar[r]^{pr_1} \ar[d] & X \\
               & \C^* & }\]

\end{proposition}
\begin{proof}
  We are working with a small enough representative of the germ $(X,0)\subset (\C^{n+1},0)$ embedded in such
a way that $Y \subset X$ is linear, this implies that we will have:
\begin{enumerate}
 \item $C(X) \subset \C^{n+1} \times \check{\P}^n$
 \item $\X \subset \C^{n+1} \times \C$.
 \item $C(\X) \subset \C^{n+1} \times \C \times \check{\P}^{n+1}$
 \item $C_\varphi(\X) \subset \C^{n+1} \times \C \times \check{\P}^{n}$
\end{enumerate}
  We will actually work with the non-projectivized versions of the conormal spaces, that is with the spaces
 $T^*_X(\C^{n+1})$, $T^*_\X(\C^{n+1}\times \C)$ and $T^*_\X((\C^{n+1}\times \C)/\C)$ respectively. Moreover, we 
  will fix a coordinate system \linebreak $(z_0, \ldots, z_{n-s},y_1,\ldots,y_s, t, a_0, \ldots,a_{n-s},c_1,\ldots,
 c_s,b)$ of $\C^{n+1} \times \C \times \check{\C}^{n+1} \times \check{\C}$. By construction,
 the map $\phi: \X \to X \times \C$ is an isomorphism when restricted to $\X \setminus \X(0)$ and
 has $X \times \C^*$ as its image. Actually, this alone implies that both the conormal space $C(\X \setminus \X(0))$ 
 and the relative conormal space $C_\varphi(\X \setminus \X(0))$ are isomorphic to $C(X) \times \C^*$. However
 to verify that we have the commutative diagram we will specify these isomorphisms. Recall that the series
 \[F_i=t^{-m_i}f_i(tz_0,\ldots,tz_{n-s},y_1,\ldots,y_s), \; i=1,\ldots,p\]
 define the specialization space $\X$ in $\C^{n+1} \times \C$, where $m_i = \nu_Y f_i$. \\

  Let $x=(z,y,t)$, $t \neq 0$, be a smooth point of $\X$, then it is a smooth point of $\X(t)$, and $\phi(x)=
(tz,y,t)$ is a smooth point of $X \times \C^*$; consequently $(tz,y)$ is a smooth point of $X$. Now, for 
any point $(x,a,c,b)$ in $\kappa_\X^{-1}(x)$ we have that there exist constants $\lambda_1, \ldots, 
\lambda_p$ such that:
\begin{align}\label{coordinatesa}
   a_j & = \sum_{i=1}^p \lambda_i \frac{\partial F_i}{\partial z_j}(x) = \sum_{i=1}^p \lambda_i t^{-m_i + 1}
          \frac{\partial f_i}{\partial z_j} (tz,y) \\
   c_j & = \sum_{i=1}^p \lambda_i \frac{\partial F_i}{\partial y_j}(x) = \sum_{i=1}^p \lambda_i t^{-m_i}
           \frac{\partial f_i}{ \partial y_j} (tz,y) \label{coordinatesc}\\
     b & = \sum_{i=1}^p \lambda_i \frac{\partial F_i}{\partial t}(x)= \sum_{i=1}^p \lambda_i \left( (-m_i)t^{-m_i+1}f_i(tz,y)
    + t^{-m_i}(\sum_{k=0}^{n-s} z_k \frac{\partial f_i}{\partial z_k}(tz,y))\right) \\
       & = \sum_{i=1}^p \lambda_i \left( t^{-m_i}(\sum_{k=0}^{n-s} z_k \frac{\partial f_i}{\partial z_k}(tz,y))\right),
    \; \; \textrm{because} \; f_i(tz,y)=0 \; \textrm{on} \; X \times \C.\label{coordinatesb}  
\end{align} 
  Analogously, for any point $(x,a,c)$ in $\kappa_\varphi^{-1}(x)$, there exist constants 
$\lambda_1, \ldots, \lambda_p$ such that, the coordinates $a_j$ and $c_j$ are given by the corresponding 
equations \ref{coordinatesa} and \ref{coordinatesc}. This implies that the natural projection $P:(z,y,t,a,c,b) 
\mapsto (z,y,t,a,c)$ induces a surjective morphism to $C_\varphi(\X \setminus \X(0))$ when restricted to
$C(\X \setminus \X(0))$. But, from \ref{coordinatesb} we can see that $tb=\sum_{k=0}^{n-s}z_k a_k$, so as long
as $t\neq 0$ the $b$ coordinate is completely determined by the $a$ and $z$ coordinates which proves that the 
aforementioned map $P$ is an isomorphism.\\

   On the other hand, for the corresponding point $x'=(tz,y)$ of $X$, we have that for any point $(x',a,c)$ in 
$\kappa_X^{-1}(x')$ there exists constants $\alpha_1, \ldots, \alpha_p$ such that:
\begin{align*}
   a_j & = \sum_{i=1}^p \alpha_i \frac{\partial f_i}{\partial z_j} (tz,y) \\
   c_j & = \sum_{i=1}^p \alpha_i \frac{\partial f_i}{ \partial y_j} (tz,y) \\
\end{align*}
  This implies that if $t \neq 0$, the automorphism 
 $\Upsilon: \C^{n+1} \times \C \times \check{\C}^{n+1} \circlearrowleft $ of the ambient space defined by: 
 \[(z,y,t,a,c) \mapsto (tz_0,\ldots,tz_{n-s},y_1,\ldots,y_s,t,a_0,\ldots,a_{n-s},tc_1,\ldots,tc_s)\]
 induces a surjective map $\psi_\varphi:C_\varphi(\X \setminus \X(0)) \to C(X) \times \C^* $ simply by 
 setting $\lambda_i = t^{m_i-1}\alpha_i$. Moreover, since the map $\Upsilon$ is biholomorphic in the 
 open dense set $t\neq 0$, the map $\psi_\varphi$ is an isomorphism. 
\end{proof}

\begin{remark}\label{extension} In regard to the previous diagrams, note that:
  \begin{enumerate}
	 \item The map $\phi$ is defined on all of $\X$, and the image of the special fiber $\X(0)$
	      is just the origin in $X \times \C$. Note as well, that for a fixed $t \neq 0$,  the 
	      morphism $pr_1 \circ \phi| : \X(t) \to X$ is an isomorphism.
	 \item The obstruction to the extension of $\psi$ to $C(\X)$ comes from the map $\check{\P}^{n+1} \to 
	       \check{\P}^n$, which is undefined at the point $[0:\cdots:0:1]$. This means that for any point
	       $\left( (\z,t), [\a:b]\right)$ in $C(\X) \cap \left(\X \times \{\check{\P}^{n+1}\setminus [0:1]\}
	       \right)$, the hyperplane $[\a] \in \check{\P}^n$ is tangent to $X$ at the point $t\z=(tz_0,
	       \ldots, tz_n)$. In particular, for $t=0$ the hyperplane $[\a]$ is tangent to $X$ at the 
	       origin. 
  \end{enumerate}   
  
 \end{remark} 

\section{The Normal/Conormal diagram.}    

    Let $(Y,0)\subset (X,0)$ be a germ of nonsingular analytic subvariety of dimension  $s < d$ as 
  before. The Whitney conditions of the pair $(X^0,Y)$ at $0$ can be expressed in terms of the normal/conormal 
  diagram of the pair $(X,Y,0)$. We will choose an embedding $(X,0) \subset (\C^{n+1},0)$ such that 
  the germ $(Y,0)$ is linear with coordinate system $(z_0, \ldots, z_{n-s}, y_1, \ldots, y_s)$. 
      
     \[\xymatrix{E_YC(X)\ar[r]^{\hat{e}_Y}\ar[dd]^{\kappa'_X}\ar[ddr]^\zeta & C(X)\ar[dd]^\kappa_X
             \\
             & &  \\
             E_YX\ar[r]_{e_Y}  &   X }\]

   We will denote by $r:(X,0)\longrightarrow (Y,0)$ the retraction induced by the projection 
  onto the $y$ coordinates.
  \newpage

 \begin{proposition}\label{WhitneyinNormalConormal}
      Let $D$ denote the reduced divisor $|\zeta^{-1}(Y)| \subset E_YC(X)$, then:
  \begin{enumerate}
    \item The pair $(X^0,Y)$ satisfies Whitney's condition a) at every 
     point $y \in Y$ if and only if we have the set theoretical equality 
    $|C(X)\cap C(Y)| = |\kappa_X^{-1}(Y)|$. 

     \item The pair $(X^0,Y)$ satisfies Whitney's condition a) at every 
     point $y \in Y$ if and only if $D$ is contained in $Y \times \P^{n-s} \times \check{\P}^{n-s}$
     where for every $y \in Y$, $\check{\P}^{n-s}$ denotes the space of hyperplanes containing $T_yY$.
     In particular, they satisfy Whitney's condition a) at $0$ if and only if $\zeta^{-1} (0) \subset
     \{0\} \times \P^{n-s}\times \check {\P}^{n-s}$.

     \item The pair $(X^0,Y)$ satisfies Whitney's condition b) at $y \in Y$ if and
  only if $|\zeta^{-1}(y)|$ is contained in the incidence variety $I \subset \{y\} \times \P^{n-s} \times 
  \check{\P}^{n-s}$. 
  \end{enumerate}
  \end{proposition}
  \begin{proof}

       Whitney conditions are defined in terms of limit of tangent spaces. However, once we have fixed an 
     embedding $(X,0) \subset (\C^{n+1},0)$, since a hyperplane $H$ is a limit of tangent hyperplanes
    if and only if it contains a limit of tangent spaces  we can restate Whitney conditions:

    \begin{itemize}
     \item[$\bullet$)] The pair $(X^0,Y)_0$ satisfies Whitney condition a) at $0$ if for any sequence 
               of non singular points $\{x_i\}_{i \in \N}\subset X^0$ tending to $0$, and any
               sequence $\{H_i\}_{i \in \N}$ where $H_i$ is a tangent hyperplane to $X$ at the point 
               $x_i$ we have the inclusion
               \[ T_0Y \subset \lim_{i\to \infty} H_i \]
     \item[$\bullet$)] The pair $(X^0,Y)$ satisfies Whitney condition b) at $y \in Y$ if for any sequence
               of non singular points $\{x_i\}_{i \in \N}\subset X^0$ tending to $y$, and any
               sequence $\{H_i\}_{_ \in \N}$ where $H_i$ is a tangent hyperplane to $X$ at the point 
               $x_i$ we have the inclusion
               \[ \lim_{i\to \infty} [x_ir(x_i)] \subset \lim_{i\to \infty} H_i\]       
\end{itemize}

       With this in mind \textbf{1)} is now only an observation. Note that we always have the inclusion 
     $|C(X)\cap C(Y)| \subset |\kappa^{-1}(Y)|$. On the other hand, the inclusion 
     $|\kappa^{-1}(Y)|\subset |C(Y)|$ means that for every $y \in Y$ every limit of tangent hyperplanes 
     to $X$ at $y$,  $H \in \kappa_X^{-1}(y)$, is also a tangent hyperplane to $Y$ at $y$, that is $T_yY 
    \subset H$.\\  
 
       For \textbf{2)}, with the coordinate system we have fixed we have 
     the blowing up $E_YX$ as a subspace of $X \times \P^{n-s}$, and the conormal space $C(Y)$  equal
     to $Y \times \check{\P}^{n-s}$ where $\check{\P}^{n-s}$ corresponds to the projective dual of $\P Y$, 
     that is the algebraic set defined by $c_1= \cdots=c_s=0$. Then, from \textbf{1)}
     satisfying condition a) is equivalent to the inclusion $|\kappa^{-1}_X(Y)| \subset
     Y \times \check{\P}^{n-s}$ which by construction of the normal conormal diagram is equivalent
     to the inclusion $|\zeta^{-1} (Y)| \subset Y \times \P^{n-s}\times \check{\P}^{n-s}$.\\

       To prove \textbf{3)}, with the coordinate system we have fixed, we have the natural 
  retraction $r:\C^{n+1} \to Y$ sending $(\z,\underline{y}) \to y$ which at the same time is used
  to build  the underlying set of the blowup of $X$ along $Y$, $E_YX$. So, from the 
  construction of $E_YC(X)$ as a subspace of the fiber product, we have to take the closure of the   
  set of points of this space of the form $(\z,\underline{y},l,H)$ where  $(\z,\underline{y})$ is a point 
  in $X^0 \setminus Y$, $l \in \P^{n-s}$ is the line defined by $[(\z,\underline{y})-r(\z,\underline{y})]$ 
  and $H$ is a tangent hyperplane to $X$ at the point $(\z,\underline{y})$. Then, a point in the divisor 
  $D = \zeta^{-1}(Y)$ is a point 	
  $(0,y,l,H)$, where $(0,y)$ is a point in $Y$, and $l$ and $H$ are a line and a hyperplane obtained 
  in the way described in the definition of condition b) above. Finally the inclusion $l \subset H$ 	
  is just what it means that the pair $(l,H)$ is in  the incidence variety $I \subset \{y\} \times \P^{n-s} 
  \times \check{\P}^{n-s}$, which finishes the proof. 
   \end{proof}
    
 \section{Whitney's conditions.}  

   Let $(X,0) \subset (\C^{n+1},0)$ be a reduced germ of analytic singularity of
  pure dimension $d$, and let $\varphi: (\X,0) \to (\C,0)$ denote the specialization of $X$ to its 
  tangent cone $C_{X,0}$. Let $\X^0$ denote the open set of smooth points of $\X$, and let $Y$ denote
  the smooth subspace $0 \times \C \subset \X$.  Our aim is to study the equisingularity of $\X$ along $Y$, that 
  is, we want to determine whether it is possible to find a  Whitney stratification of $\X$  
  in which the t-axis $Y$ is a stratum.\\

  The first step to find out if such a stratification is possible, is to verify that the pair 
  $(\X^0, Y)$ satisfies Whitney's conditions. Since $\X \setminus \X(0)$ is isomorphic to the product 
  $X \times \C^*$, Whitney's conditions are automatically verified everywhere in $\{0\} \times \C$ , with 
  the possible exception of the origin. The following result tells us that in this particular case
  it is enough to check for Whitney's condition a).  
  
  \begin{proposition}\label{aimplicab}
     If the pair $(\X^0,Y)$ satisfies Whitney's condition a) at the origin, then
   it also satisfies Whitney's condition b) at the origin.  
  \end{proposition}

   Before proving proposition \ref{aimplicab}, we need the following lemma.
  
      \begin{lema}\label{NorTan}
         There exists a natural morphism $\omega: E_Y\X \to E_0X$, making the following diagram
         commute: 
      \[\xymatrix{ E_Y\X \ar[r]^\omega \ar[d]_{e_Y} &  E_0X \ar[d]^{e_o}\\
                   \X  \ar[r]_\phi & X }\]
        Moreover, when restricted to the exceptional divisor $e_Y^{-1}(Y)= \P C_{\X,Y}$ it induces
        the natural map $\P C_{\X,Y}= Y \times \P C_{X,0} \to \P C_{X,0}$.
      \end{lema}
      \begin{proof}
          Algebraically, this results from the universal property of the blowup $E_0X$. We start with 			the diagram:
           \[\xymatrix{ E_Y\X  \ar[d]_{e_Y} &  E_0X \ar[d]^{e_o}\\
                   \X  \ar[r]_\phi & X }\]
      In this coordinate system, the maximal ideal $\m$ of the analytic algebra $O_{X,0}$ 
      is generated by $\left<z_0, \ldots, z_n\right>$. The map $\phi$, induces a morphism
      of analytic algebras $O_{X,0} \to O_{\X,0}$ defined by $z_i \mapsto tz_i$. So we have to
      prove that the ideal $\left< tz_0, \ldots, tz_n \right> \subset O_{\X,0}$ is locally invertible
      when pulled back to $E_Y\X$. But as ideals we have the equality $\left<tz_0, \ldots,tz_n 
      \right>= \left< t \right> \centerdot \left< z_0, \ldots z_n \right>$. And by definition of the
      blowup, the ideal $\left< z_0, \ldots, z_n\right > \subset O_{\X,0}$ corresponding to $Y$ is
      locally invertible when pulled back 
      to $E_Y\X$. After multiplication by a invertible ideal, it will remain 
      locally invertible. Note that, for the diagram to be commutative the morphism $\omega$ must map
      the point $(z,t), [z] \in E_Y\X \setminus \{Y \times \P^n\} \subset \X \times \P^n$ 
      to the point $(tz), [z] \in E_0X \subset X \times \P^n $ and the result follows.       
      \end{proof}    

     \begin{remark}\label{puntosclaves} Note that: 
     \begin{enumerate}
     \item For any point $y \in Y$, the tangent cone $C_{\X,y}$ is isomorphic to $C_{X,0} \times Y$,
     and the isomorphism is uniquely determined once we have chosen a set of coordinates. The 
     reason is that for any $f(z)$ vanishing on $(X,0)$, the function 
      $F(z,t)=t^{-m}f(tz)=f_m+tf_{m+1}+ t^2f_{m+2}+ \ldots$, vanishes in $(\X,0)$ and so for any
      point $y=(\underline{0},t_0)$ the initial form of $F(z,t+t_0)$ in $\C\{z_0, \ldots, z_n,
      t\}$ is equal to the initial form of $f$ at $0$. That is in$_{(0,t_0)}F=\textrm{in}_0f$.
      
     \item The projectivized normal cone $\P C_{\X,Y}$ is isomorphic to $Y \times \P C_{X,0}$. 
            This can be seen from the equations used to define $\X$ (Chapter 1, eq. \ref{eq1}), 
            where the initial form of $F_i$ with respect to $Y$, is equal to the initial form of $f_i$ at 
            the origin. That is in$_YF_i=\textrm{in}_0f_i$.
     
   \end{enumerate}
 \end{remark} 
     
     Now we can proceed to the proof of \ref{aimplicab}.
    
    \begin{proof} (\textit{ Proposition \ref{aimplicab}})
	
       We want to prove that the pair $(\X^0,Y)$ satisfies Whitney's condition b) at the origin. 
    We are assuming that it already satisfies condition a), so in particular we have that 
    $\zeta^{-1}(0)$ is contained in $\{0\} \times \P^n \times \check{\P}^n$. By proposition \ref{WhitneyinNormalConormal} 
    it suffices to prove that any point $(0,l,H) \in \zeta^{-1}(0)$  is contained in the incidence variety 
    $I \subset \{0\} \times \P^n \times \check{\P}^n$. Consider the diagram:

   \[\xymatrix{E_YC(\X)\ar[r]^{\hat{e}_Y}\ar[dd]^{\kappa'_\X}\ar[ddr]^\zeta & C(\X)\ar[dd]^{\kappa_\X}\ar[r]^{\psi} & 
             C(X) \times \C \\
             & &  \\
             E_Y\X\ar[r]_{e_Y} \ar[d]_{\omega}  &  \X \\
             E_0X & & }\]
     
       By construction, there is a sequence $(z_m, t_m, l_m, H_m)$ in $E_YC(\X) \hookrightarrow
   C(\X) \times_\X E_Y\X $ tending to $(0,l,H)$ where $(z_m,t_m)$ is not in $Y$. Through $\kappa'_\X$, we obtain a 
   sequence $(z_m,t_m,l_m)$ in $E_Y\X$ tending to $(0,l)$, and through $\hat{e}_Y$ a sequence $(z_m,t_m,H_m)$
   tending to $(0,H)$ in $C(\X)$. \\

       Now,  using the notation of proposition \ref{IsomorfismosConormales}, through the map $\psi$ 
   we obtain the sequence $(t_mz_m, \widetilde{H}_m)$ and since by hypothesis we have $b=0$, then 
   by remark \ref{extension}-2 both the sequence and its limit $(0,\widetilde{H})$ are in $C(X)$. 
   Note that if $H$ has coordinates $[a_0: \cdots: a_n:0]$, then $\widetilde{H}= 
   [a_0: \cdots: a_n] \in \check{\P}^n$. 
       On the other hand, by lemma \ref{NorTan} we have that both the
   sequence $(t_mz_m, l_m)$ obtained through the map $\omega$ and its limit $(0,l)$ are in $E_0X$. 
   Finally, Whitney's lemma (\cite[Thm. 22.1, p. 547]{Whi1} or \cite[Thm. 1.1.1]{Le1}) tells us that in 
   this situation we have that $l \subset \widetilde{H}$ and so the point $(0,l,H)$ is in the incidence variety.
 
     If the sequence $(z_m, t_m, l_m, H_m)$ in $E_YC(\X)$ is contained in the special fiber, that is $t_m=0$ 
  for all $m$, then either the point $(z_m, 0)$ is a smooth point of $\X$ and so the line $l_m=[z_m:0]$ is 
  contained in every tangent hyperplane $H_m$, or it is a singular point of $\X$ and by constructing a 
  sequence of smooth points in $\X \setminus \X(0)$ tending to it and using the maps $\psi$ and $\omega$ 
  as we did before we prove that the line $l_m$ is contained in $H_m$. In any case, what we have is that for any 
  point in the sequence $(z_m, 0, l_m, H_m)$ we already have the inclusion $l_m \subset H_m$ and so the limit
  $(0,l,H)$ satisfies this condition as well.
   \end{proof}

      The following result tells us that in order to have $Y$ be a stratum in a Whitney stratification 
    of $\X$, the condition of $(X,0)$ not having exceptional cones is necessary. \\

    \begin{lema}\label{GeneralizationNecessity}
     Let $(X,0) \subset (\C^{n+1},0)$ be a reduced germ of analytic singularity of
  pure dimension $d$, and let $\varphi: (\X,0) \to (\C,0)$ denote the specialization of $X$ to its 
  tangent cone $C_{X,0}$. Let $\X^0$ denote the open set of smooth points of $\X$, and let $Y$ denote
  the smooth subspace $0 \times \C \subset \X$. If the tangent cone $C_{X,0}$ is reduced and the pair $(\X^0, Y)$
  satisfies Whitney's condition a) then the germ $(X,0)$ does not have exceptional cones.
    \end{lema}
   \begin{proof}
      First of all, by hypothesis the pair $(\X^0, Y)$ satisfies Whitney's condition a), so by proposition
      \ref{aimplicab} it also satisfies Whitney's condition b).
      Recall that the aureole of $(\X,0)$ along $Y$ is a collection $\{V_\alpha\}$ of subcones of the normal 
      cone $C_{\X,Y}$ whose projective duals determine the set of limits of tangent hyperplanes to $\X$ at
      the points of $Y$ in the case that the pair $(\X^0,Y)$ satisfies Whitney conditions a) and b) at every point
      of $Y$ (See \cite[Thm. 2.1.1, Coro 2.1.2 p. 559-561]{L-T2} ). 
      Among the $V_\alpha$ there are the irreducible components of $|C_{\X,Y}|$. Moreover:
      \begin{enumerate}
	       \item By remark \ref{puntosclaves} we have that $C_{\X,Y} = Y \times C_{X,0}$ so its 
	        irreducible components are of the form $ Y \times \widetilde{V}_\beta$ where 
	        $\widetilde{V}_\beta$ is an irreducible component of $|C_{X,0}|$.
	        \item For each $\alpha$ the projection $V_\alpha \to Y$ is surjective and all the fibers 
	         are of the same dimension. (See \cite{L-T2}[Proposition 2.2.4.2, p. 570])
	        \item The hyperplane $H=[0:0: \cdots:1] \in \check{\P}^{n+1}$ is transversal to $(\X,0)$ by
	         hypothesis, and so by \cite[Thm. 2.3.2, p. 572]{L-T2} the collection $\{V_\alpha \cap H\}$
	         is the aureole of $\X \cap H$ along $Y \cap H$. 
	     \end{enumerate}
        
      Notice that $(\X \cap H,Y\cap H)$ is equal to $(\X(0),0)$, which is isomorphic to the tangent cone 
      $(C_{X,0},0)$ and therefore does not have exceptional cones. This means that for each $\alpha$ either
      $V_\alpha \cap H$ is an irreducible component of $C_{X,0}$ or it is empty. But the intersection can't be 
      empty because the projections $V_\alpha \to Y$ are surjective. Finally since all the fibers of the
      projection are of the same dimension then the $V_\alpha$'s are only the irreducible components of 
      $C_{\X,Y}$. But this means, that if we define the affine hyperplane $H_t$ as the hyperplane 
      with the same direction as $H$ and passing through the point $y=(0,t) \in Y$ for $t$ small enough; 
      $H_t$ is transversal to $(\X,y)$ and so we have again that the collection $\{V_\alpha \cap H_t\}$
      is the aureole of $\X \cap H_t$ along $Y \cap H_t$, that is the aureole of $(X,0)$, so it does not have
      exceptional cones.
    
   \end{proof}

      We can now use lemma \ref{GeneralizationNecessity} to prove that the Whitney conditions of the pair 
    $(\X^0,Y)$ imply that the germ $(\X,0)$ does not have exceptional cones.
  
    \begin{proposition}\label{PropoGeneral*} 
    Let $(X,0) \subset (\C^{n+1},0)$ be a reduced germ of analytic singularity of
  pure dimension $d$, and let $\varphi: (\X,0) \to (\C,0)$ denote the specialization of $X$ to its 
  tangent cone $C_{X,0}$. Let $\X^0$ denote the open set of smooth points of $\X$, and let $Y$ denote
  the smooth subspace $0 \times \C \subset \X$. 
\begin{enumerate}
	\item If the germ $(\X,0)$ does not have exceptional cones, then the pair $(\X^0, Y)$ satisfies 
	Whitney's condition a) at the origin.
	\item Moreover, if the tangent cone $C_{X,0}$ is reduced and the pair $(\X^0, Y)$
	 satisfies Whitney's condition a) at the origin then $(\X,0)$ does not have exceptional cones.
\end{enumerate}
\end{proposition}

\begin{proof}
Let us choose a representative of $(X,0)$ in $(\C^{n+1},0)$, then $(\X,0) \subset (\C^{n+2},0)$.
Let $C(\X) \subset \C^{n+2} \times \check{\P}^{n+1}$ denote the conormal space of $\X$, and let us consider the following diagram:
   
  \[\xymatrix { C(\X) \ar[d]^{\kappa_\X}  & C(Y)\ar[d]^h
                   \\
                  \X &  Y\ar @{_{(}->}[l]  }\]
   By proposition \ref{WhitneyinNormalConormal}, Whitney's condition a) at the origin is equivalent to 
   the set theoretic inclusion 
   \[|\kappa_\X^{-1}(0)| \subset |h^{-1}(0)|\]
   Let $\left( (z_0,\ldots,z_n,t), [a_0:a_1:\ldots:a_n:b]\right)$ be the coordinates of 
   $\C^{n+2} \times \check{\P}^{n+1}$ as before. Now, since $Y$ is the $t$ axis, the conormal 
   space $C(Y)$ is defined by the equations $z_0=\cdots=z_n=b=0$, and for $h^{-1}(0)$ we just add 
   the equation 
   $t=0$. \\
   
   1) By hypothesis $(\X,0)$ does not have exceptional cones, which means that $|\kappa_\X^{-1}(0)|$ is 
      just the dual of the tangent cone $C_{\X,0}= C_{X,0} \times \C$. In particular, every tangent 
      hyperplane to $C_{\X,0}$ contains the $t$ axis, that is $b=0$, so is contained in $h^{-1}(0)$,
      and we have Whitney's condition a).\\
      
   2 ) By lemma \ref{GeneralizationNecessity} we know that $(X,0)$ does not have exceptional cones. 
      Since every point in
      $\kappa_{\X}^{-1}(0)$, that is every tangent hyperplane to
      $\X$ at the origin satisfies $b=0$, the remark \ref{extension}-2 tells us that the
      morphism $(\widetilde{pr_1} \circ \psi):C(\X \setminus \X(0)) \to C(X)$ of proposition 
      \ref{IsomorfismosConormales}, sending $(z,t),[a:b] \to (tz),[a]$ can be extended 
      to $C(\X)$.
      In particular the point, $(0), [a]$ is in $\kappa_X^{-1}(0) \subset C(X)$, and since $(X,0)$ does not have
      exceptional
      cones, then $[a]$ is in the dual of the tangent cone $C_{X,0}$, which implies that 
      $\kappa_{\X}^{-1}(0)$ is just the dual of the tangent cone $C_{\X,0}$, and $(\X,0)$
      does not have exceptional cones. 
           
   \end{proof}    

     We will study Whitney's condition a) by deriving a characterization specific to our situation 
   from the characterization given first by Teissier in \cite{Te2} in the case of isolated hypersurface 
   singularities and subsequently generalized by Gaffney in \cite{Ga2} in terms of integral dependence 
   of modules.\\

\section{Limits of tangents spaces and integral closure of modules}

    There are several equivalent definitions of integral closure for modules. In our case, it is simpler
to work with the following definition, as stated in \cite[Section 3, p. 555]{GaKl}. 

\begin{definition}\label{DefIntClosMod}
   Let $O_{\X}^p$ be a free module of rank $p \geq 1$. Let $M$ be a coherent submodule of 
$O_{\X}^p$ and $h \in O_{\X}^p$. Given a map of germs $\phi: (\C,0) \to (\X,0)$, denote by $h \circ 
\phi$ the induced section of $O_1^p$, and by $M \circ \phi$ the induced submodule. Call $h$ \textit{ 
integrally dependent (resp. strictly dependent) on $M$} at $0$ if, for every $\phi$, the section $h \circ 
\phi \in O_1^p$ belongs to the submodule $M \circ \phi$ of $O_1^p$ (resp. to the submodule $\m_1(M \circ
\phi)$), where $\m_1$ is the maximal ideal of $O_1 = \C\{\tau\}$. The submodule of $O_{\X}^p$
generated by all such $h$ will be denoted by $\overline{M}$, respectively by $M^\dag$. 

  Moreover, we say that a submodule $N \subset M$ is a \textbf{reduction} of $M$ if $\overline{N}=\overline{M}$. 
\end{definition}
   
    If the germ $(\X,0)$ is not irreducible, for every irreducible component $\X_i$ of $\X$ the module $M$ induces 
  a submodule $M_{\X_i}$ of $O_{\X_i}^p$ via the morphism of analytic algebras $O_{\X,0} \to O_{\X_i,0}$, and 
  the same goes for a section $h$ of $O_\X^p$. A simple calculation then shows:

\begin{lema}\label{DependenceIrreducible}
   Let $(\X,0)= \bigcup_{i=1}^r (\X_i,0)$ be the irreducible decomposition of the germ. The section $h$ is integrally 
 dependent (respectively strictly dependent) on $M$ at $0$ if and only if for every irreducible component $\X_i$ 
 the induced section $h_i$ is integrally dependent (respectively strictly dependent) on $M_{\X_i}$ at $0$.
\end{lema}

  We will state the main results we will be using.
 Let $M$ be a coherent submodule of $O_{\X}^p$ as before, and let $[M]$ be a matrix of generators of $M$ 
for a small enough neighborhood of the origin in $(\X,0)$, that is the matrix describing the morphism $\mu$ of:
\[O_\X^q \stackrel{\mu}{\longrightarrow} O_\X^r \longrightarrow O_\X^p/M \longrightarrow 0\]
Let $J_k(M)$ denote the ideal of $O_\X$ generated by the $k \times k$ minors of $[M]$. This is the same as the 
$(p-k)$-th Fitting ideal of $O_{\X}^p/M$ and so is independent of the choice of generators of $M$. If $h \in O_{\X}^p$, let 
$(h,M)$  denote the submodule of $O_{\X}^p$ generated by $h$ and $M$. 

\begin{proposition}\label{Moduloseidealjacobiano}\cite[Prop 1.7, p. 304]{Ga1},and \cite[Prop 1.5, p. 57]{Ga2}\\
  Suppose $M$ is a submodule of $O_{\X}^p$, $h \in O_{\X}^p$ and the rank of $(h,M)$ is $k$ on each irreducible component
of $(\X,0)$. Then $h$ is integrally dependent (resp. strictly dependent) on $M$ at $0$ if and only if each minor
in $J_k(h,M)$ which depends on $h$ is integrally dependent (resp. strictly dependent) on $J_k(M)$. 
\end{proposition}

\begin{lema}\label{IdealIpsi}\cite[Lemma 3.3, p. 557]{GaKl}
   For a section $h$ of $O_{\X}^p$ to be integrally dependent, respectively strictly dependent, on $M$ at $0$,
it is necessary that for all maps:
\begin{align*}
            \phi: (\C,0) &\to (\X,0)\\
            \psi: (\C,0) & \to (\mathrm{Hom}(\C^p, \C), \lambda), \;\; \lambda\neq 0
\end{align*}
the function $\psi(h \circ \phi)$ on $\C$ belongs to the ideal $I_\psi(M \circ \phi)$ generated by applying
$\psi(\tau)$ to the generators of $M \circ \phi $, respectively to the ideal $\m_1 I_\psi(M \circ \phi)$. 

  Conversely it is sufficient that this condition is satisfied for every \break 
$\phi: (\C , \C \setminus \{0\},0) \to (\X, \X \setminus W, 0)$, where $(W,0) \subset (\X,0)$ is a proper analytic subset of $\X$.
\end{lema}

\begin{corollary}\label{CoroCrecimiento}( \cite[Proposition 1.11, p. 306]{Ga1})
   The section $h$ is integrally dependent on $M$ at $0$ if and only if for each choice of generators
 $\{m_i\}$ of $M$ there exists a neighborhood $U$ of $0$ in $\X$, and a real constant $C$, such that for every 
 section $\Psi: \X \to \X \times \check{\P}^{p-1}$ of the trivial bundle $\X \times \check{\P}^{p-1}$
 and every point $z \in U$ we have:
 \[\left| \Psi(z) \cdot h(z)\right| \leq C \sup_i\left| \Psi(z) \cdot m_i(z)\right|\] 
\end{corollary}

 The previous results direct us to work with the space $\X \times \check{\C}^p$, or even with
the space $\X \times \check{\P}^{p-1}$ since we ask that the image of $\psi$ does 
not contain the point $0$ in $\check{\C}^p$. These spaces can be seen respectively as the analytic spectrum (analytic proj) of
the symmetric algebra of $O_\X^p$, that is $O_\X [u_1,\ldots, u_p]$. The section $h \in O_\X^p$ and the 
submodule $M \subset O_\X^p$ generate ideals in $O_\X[u_1, \ldots, u_p]$ which we will denote by 
$\rho(h)$ and $\rho(M)$.

\begin{remark}\label{Encajegrado1}
 Recall that the embedding of $O_\X^p$ in $O_\X[u_1, \ldots, u_p]$ is in degree 1, and is given by 
\[h = \begin{pmatrix} h_1 \\ h_2 \\ \vdots \\ h_p \end{pmatrix} \mapsto \rho(h)= u_1h_1 + \cdots + u_ph_p \]
\end{remark}

We will consider the normalized  blowup of $\X \times \check{\P}^{p-1}$ along the subspace $Z$ 
defined by the ideal $\rho(M)O_\X[u_1,\ldots, u_p]$ which we will denote by
\[\pi: \overline{E_Z(\X \times \check{\P}^{p-1})} \to \X \times \check{\P}^{p-1} \to \X\]
Its exceptional divisor will be denoted by $F$

\begin{proposition}\label{DependenciaMeromorfa}\cite[Prop. 3.5, p. 558]{GaKl} Let $h \in O_{\X}^p$, and let $Y$ be a closed 
analytic subset of the image of $F$ in $\X$. Then:
\begin{enumerate}
 \item $h$ is integrally dependent on $M$ at $0$ if and only if along each irreducible component of $F$, the ideal
       $\rho(h) \circ \pi$ vanishes to order at least the order of vanishing of $\rho(M) \circ \pi$.
 \item $h$ is strictly dependent on $M$ at every $y \in Y$ if and only if along each component $V$ of
       $F$, the ideal $\rho(h) \circ \pi$ lies in the product $I(Y,V)\rho(M)\circ \pi$, where $I(Y,V)$
       denotes the ideal of the reduced preimage of $Y$ in $V$.
\end{enumerate}
\end{proposition}

From this point on we will assume that \textbf{the germ $\mathbf{(X,0)}$ is irreducible.}\\
 
   Let $\left< F_1, \ldots, F_p\right>\C\{z_0, \ldots, z_n, t\}$ be the ideal defining the germ $(\X,0)$ as before. 
  In other words, $(\X,0) = (F^{-1}(0),0)$ where $F=(F_1, \ldots, F_p): (\C^{n+1} \times \C,0) \to (\C^p,0)$. Let 
  $c$ denote the codimension of $\X$ in $\C^{n+1} \times \C$, and let $S$ denote the set 
  of increasing sequences of $c$ positive integers less than $p+1$. For $\alpha \in S$ denote
  by $[DF]_\alpha$ the $c \times (n+2)$ submatrix of $[DF]$ formed by the $(\alpha_1,\ldots,\alpha_c)$
  lines of $[DF]$. That is the jacobian matrix, of the map $F_\alpha:=(F_{\alpha_1}, \ldots F_{\alpha_c}):
  \C^{n+1} \times \C \to \C^c$. 
 
  \begin{definition}\label{alphaDefJacMod}
   For  $\alpha \in S$, define the \textbf{$\mathbf{\alpha}$-Jacobian module of $F$} as the submodule 
$JM(F)_\alpha$ of $O_\X^c$ generated by the columns of the matrix $[DF]_\alpha$, that is:
\[ JM(F)_\alpha:= O_\X  \begin{pmatrix} \frac{\partial F_{\alpha_1}}{\partial z_0} \\ \vdots \\  
              \frac{\partial F_{\alpha_c}}{\partial z_0}  \end{pmatrix} + \cdots + 
   O_\X  \begin{pmatrix}  \frac{\partial F_{\alpha_1}}{\partial z_n} \\ \vdots \\  
             \frac{\partial F_{\alpha_c}}{\partial z_n} \end{pmatrix} + 
   O_\X  \begin{pmatrix}  \frac{\partial F_{\alpha_1}}{\partial t}  \\ \vdots \\  
           \frac{\partial F_{\alpha_c}}{\partial t}  \end{pmatrix} \subset O_\X^c  \] 
\end{definition}

 Let $v$  be a vector in $\C^{n+1} \times \C$, then by $\frac{\partial F_\alpha}{\partial v}$ we mean
the directional derivative of $F_\alpha$ with respect to $v$. That is:
\[ \frac{\partial F_\alpha}{\partial v}:= \left[DF \right]_\alpha(v)\]
 In particular $\frac{\partial F_\alpha}{\partial v}$ is a linear combination of the columns of 
 $F_\alpha$ and so it belongs to the $\alpha$-jacobian module $JM(F)_\alpha$.
 
\begin{definition}\label{Reljacmoddef}
    Given an analytic map germ $g:(\C^{n+1} \times \C, 0) \to (\C^l,0)$, and $\alpha\in S$, 
let $JM_g(F)_\alpha$ denote the submodule of $JM(F)_\alpha$ generated by the "partials" $\frac{\partial F_\alpha}{\partial v}$ for all vector fields 
$v$ on $\C^{n+1} \times \C$ tangent to the fibers of $g$, that is, for all $v$ that map to the $0-$field on $\C^l$.
Call $JM_g(F)_\alpha$ the $\mathbf{\alpha}-$\textbf{Relative Jacobian Module} with respect to $g$. 
\end{definition}

  Note that if $H$ is a hyperplane in $\C^{n+1} \times \C$ defined by the kernel of the 
linear map $h:\C^{n+1} \times \C \to \C$,
then $JM_h(F)_\alpha$ is the submodule of $JM(F)_\alpha$ generated by the partials 
$\frac{\partial F_\alpha}{ \partial v}$ for all vectors $v \in H$.\\

\begin{remark}\label{ChoosingAlpha}
  \begin{enumerate}
   \item For every non singular point $(z,t) \in \X^0$, the matrix $[DF(z,t)]$ has rank $c:=n+1-d$, and so there exists
 an $\alpha \in S$ such that at least one of the maximal minors ($c \times c$) of the matrix $[DF]_\alpha$
 is not identically zero in $O_{\X,0}$. 
   \item For every point $(z,t)$ in the relative smooth locus $\X^0_\varphi:=\bigcup\X(t)^0$, the matrix 
   $[D_\varphi F(z,t)]= \left[\frac{\partial F_i}{\partial z_j}\right]^{i=1\ldots p}_{j=0 \ldots n}$ has 
  rank $c:=n+1-d$, and so there exists a $\gamma \in S$ such that at least one of the maximal minors ($c \times c$)
  of the matrix $[D_\varphi F]_\gamma$ is not identically zero in $O_{\X,0}$.
  \end{enumerate}
\end{remark}
   
   The $\alpha-$relative jacobian module, for an appropriately chosen $\alpha$, can be used to study 
the limits of tangent hyperplanes.

     \begin{proposition}\label{AlphaExplosionZ}
   Let $\alpha \in S$ be as in remark \ref{ChoosingAlpha}-1, and let 
$E_Z(\X \times \check{\P}^{c-1}) \subset \X \times \check{\P}^{c-1} \times \check{\P}^{n+1}$ be the blowup
of $\X \times \check{\P}^{c-1}$ along the subspace $Z$ defined by the ideal $\rho(JM(F)_\alpha)O_\X[u_1,\ldots, u_c]$.
Then, there exists a surjective map \linebreak $\eta:E_Z(\X \times \check{\P}^{c-1}) \to C(\X)$, making 
the following diagram commutative:
\[\xymatrix { E_Z(\X \times \check{\P}^{c-1}) \ar[d]_{e_Z} \ar[r]^\eta & C(\X)  \ar[d]_{\kappa_\X} \\
                   \X \times \check{\P}^{c-1}  \ar[r] &  \X \ar [d]_\varphi  \\
                                & \C} \]
\end{proposition}
\begin{proof}
   Let $\alpha \in S$ be as in remark \ref{ChoosingAlpha}. Since $\X$ is irreducible, there exists 
 an open dense set $\X^0_\alpha \subset \X^0$, where for any point $(z,t) \in \X^0_\alpha$  
 the tangent space $T_{(z,t)}\X$ is the kernel of the matrix $[DF]_\alpha$, that is, it is obtained as 
 the intersection of the $c:=n+1-d$ hyperplanes $[\overrightarrow{dF_{\alpha_j}}(z,t)]$. Moreover, 
 since $c$ is the codimension of $\X$, any linear equation defining the tangent hyperplane 
 $H=[a:b]$ to $\X$ at $(z,t)$ is expressed as a \textbf{unique} linear combination of these $c$ hyperplanes 
 $H=[\sum \beta_j \overrightarrow{dF_{\alpha_j}}(z,t)]$,
 that is, they form a base of the fiber $\kappa_\X^{-1}(z,t)$ over $(z,t)$ in the conormal space $C(\X)$.
 So for any point $(z,t,u) \in \X \times \check{\C}^c$ 
 with $(z,t) \in \X^0_\alpha$ we have the map
\[(z,t,u) \in \X \times \check{\C}^c \mapsto  (z,t), \left[\sum_{i=1}^c u_i \overrightarrow{dF_{\alpha_i}}
   (z,t)\right]
 \in C(\X) \subset \X \times \check{\P}^{n+1}\]
 Note that this map is invariant with respect to the homotheties of $\check{\C}^c$, so it defines a map
 $\X \times \check{\P}^{c-1} \to \X \times \check{\P}^{n+1}$.

   On the other hand, from definition \ref{alphaDefJacMod} and remark \ref{Encajegrado1}, we get that the ideal
$\rho(JM(F)_\alpha)$ has the following system of homogeneous generators:
\[\rho(JM(F)_\alpha)= \left< u_1 \frac{\partial F_{\alpha_1}}{\partial z_0} + \cdots + u_c \frac{\partial F_{\alpha_c}}{\partial z_0},
         \ldots,u_1 \frac{\partial F_{\alpha_1}}{\partial t} + \cdots + u_c \frac{\partial F_{\alpha_c}}{\partial t} \right>
       O_\X[u_1,\ldots, u_c]  \]
 and so a point $(z,t,[u]) \in \X \times \check{\P}^{c-1}$  is in $Z$ if and only if
\[u_1\overrightarrow{dF_{\alpha_1}}(z,t)+ \cdots+ u_c\overrightarrow{dF_{\alpha_c}}(z,t)= \overrightarrow{0}\]
that is, $Z$ is the set of points where the previously stated map
\[(z,t,[u]) \in \X \times \check{\P}^{c-1} \mapsto  (z,t), \left[\sum_{i=1}^c u_i \overrightarrow{dF_{\alpha_i}}(z,t)\right]
 \in C(\X) \subset \X \times \check{\P}^{n+1}\]
is not defined. Thus, by blowing up the space $Z$ in this set of coordinates, we obtain the space
$E_Z(\X \times \check{\P}^{c-1}) \subset \X \times \check{\P}^{c-1} \times \check{\P}^{n+1}$ upon which
the morphism \[\eta: E_Z(\X \times \check{\P}^{c-1}) \to C(\X)\] is defined by the restriction 
to $ E_Z(\X \times \check{\P}^{c-1})$ of the projection $\X \times \check{\P}^{c-1} \times \check{\P}^{n+1} \to
\X \times \check{\P}^{n+1}$. Moreover, since for any point $(z,t) \in \X^0_\alpha$ and tangent hyperplane
$H \in \kappa_\X^{-1}(z,t)$ there exists a unique $[u] \in \check{\P}^{c-1}$ such that the point $(z,t,[u]) 
\notin Z$ and the point $(z,t,[u],H) \in E_Z(\X \times \check{\P}^{c-1})$, then the morphism $\eta$ is surjective.
\end{proof}
 
The proof of this proposition has the following result as an immediate coro- llary.

\begin{corollary}\label{CoroAlphaExplosionZ}
  For each appropriately chosen $\alpha \in S$, the restriction of 
 $\eta$ to $e_Z^{-1}(\X^0_\alpha)$ is an isomorphism. In other words, the analytic spaces
 $\X_\alpha^0 \times \check{\P}^{c-1}$ and $\kappa_\X^{-1}(\X_\alpha^0)$ are isomorphic.
\end{corollary}

\begin{remark}\label{AlphaNash}
  In the same spirit of the proof of the previous proposition, we can see that by choosing a $\gamma \in S$ as in 
remark \ref{ChoosingAlpha}-2, the irreducibility of $\X$ together with the constructive proof of 
\ref{NashNobileCoro} implies that the blowup of the ideal $J_c(JM_\varphi(F)_\gamma)$ generated
by the maximal minors of $[D_\varphi F]_\gamma$ gives the relative Nash modification $\Na_\varphi \X$.
\end{remark}

  The link between limits of tangent hyperplanes and the integral closure is further explained in the
 following results.

\begin{lema}\label{Limittangenthyperplane}(\cite[Lemma 2.1, p. 58]{Ga2})
 Let $(\X,0) \subset (\C^{n+1}\times \C,0)$ be defined by $F^{-1}(0)$ as before, let $\alpha \in S$ be as in 
remark \ref{ChoosingAlpha}-1, and let $\X_\alpha^0$ be the open dense set of smooth points of $\X$ where  
the kernel of the matrix $[DF]_\alpha$ defines the tangent space $T_{(z,t)}\X$. 
Then a  hyperplane $H=[a_0: \cdots: a_n:b] \in \check{\P}^{n+1}$ is a limit of tangent hyperplanes to $(\X,0)$ 
if and only if there exists a pair of maps $\phi:(\C,\C\setminus{0},0) \to (\X,\X_\alpha^0,0)$ and 
$\psi: (\C,0) \to (\check{\C}^c, \lambda \neq 0)$ such that the point $(\phi(\tau), \psi(\tau)) \notin Z \subset \X \times \check{\C}^c$ and
for some k
\[(a_0, \ldots , a_n, b) = \lim_{\tau \to 0} \frac{\psi(\tau)DF_\alpha(\phi(\tau))}{\tau ^k}\] 
\end{lema}
\begin{proof}
   The proof of this result is basically the same as that given in the reference just by noting the 
equivalence between a map 
\[\Theta: (\C,\C\setminus\{0\},0) \to (\X \times \C^c, \X_\alpha^0 \times \C^c, (0,\lambda))\] 
and the pair of maps $\phi:(\C,\C\setminus\{0\},0) \to (\X, \X_\alpha^0,0)$ and $\psi: (\C,0) \to (\check{\C}^c, \lambda \neq 0)$, and 
then using propositon \ref{AlphaExplosionZ} and its corollary.
\end{proof}

\begin{corollary} \label{NoHayHiperplanoVertical} Let $\varphi:(\X,0) \to \C$  denote the 
specialization of $(X,0)$ to its tangent cone $(C_{X,0},0)$. The hyperplane $\{t=0\}$ is not a limit 
of tangent hyperplanes to $\X$ at $(z,t)$ if and only if $\frac{\partial F}{\partial t} \in 
\overline{JM_\varphi (F)}$ in $O_{\X,(z,t)}$.
\end{corollary}
\begin{proof}
  From lemma \ref{Limittangenthyperplane}, the hyperplane $\{t=0\}$ is a limit of tangent
 hyperplanes if and only if there exists a pair of maps $\phi:(\C,0) \to (\X^0,(z,t))$ and 
 $\psi: (\C,0) \to (\check{\C}^p, \lambda \neq 0)$ such that the point 
 $(\phi(\tau), \psi(\tau)) \notin Z \subset \X \times \check{\C}^p$ and for some k
\[ (0, \ldots , 0, \alpha) = \lim_{\tau \to 0} \frac{\psi(\tau)DF(\phi(\tau))}{\tau ^k}\] 
But we can see that $\psi(\tau)DF(\phi(\tau))$ is equal to 
\[\left( \rho \left(\frac{\partial F}{\partial z_0}\right)(\phi(\tau),\psi(\tau)),\ldots,
    \rho \left(\frac{\partial F}{\partial z_n}\right)(\phi(\tau),\psi(\tau)),
    \rho \left(\frac{\partial F}{\partial t}\right)(\phi(\tau),\psi(\tau))  \right) \]
 and so, if we denote by $\textrm{ord}_0 \gamma(\tau)$ the order of the series $\gamma(\tau)$ in $\C\{\tau\}$, the 
 limit condition tells us that 
 \[\textrm{ord}_0 \rho \left(\frac{\partial F}{\partial t}\right)(\phi(\tau),\psi(\tau)) < 
   \textrm{ord}_0 \rho \left(\frac{\partial F}{\partial z_j}\right)(\phi(\tau),\psi(\tau)), \; 
  \textrm{for} \; j=0,\ldots,n\]
 This implies that for every $C \in \R$ there exists an $\epsilon \in \R$ such that for every
  $|\tau| < \epsilon$ we have that $|\rho \left(\frac{\partial F}{\partial t}\right)(\phi(\tau),\psi(\tau))| >
  C|\rho \left(\frac{\partial F}{\partial z_j}\right)(\phi(\tau),\psi(\tau))|$. Corollary \ref{CoroCrecimiento}
  finishes the proof.
\end{proof}

  The equivalence statement of Whitney's condition a) in terms of integral closure and the jacobian
module given by Gaffney and Kleiman(\cite[Cor 2.4, p. 60 ]{Ga2} or \cite[lemma 4.1, p. 560]{GaKl})   
can now be refined in the irreducible case by using the $\alpha-$relative jacobian module with 
basically the same proof.
  
\begin{Theorem}\label{Whitagen}$\;$\\
  Let $(X,0) \subset (\C^{n+1},0)$ be an irreducible and reduced germ of analytic singularity defined 
by an holomorphic map $f:(\C^{n+1},0) \to (\C^p,0)$, $X = f^{-1}(0)$. Let $(V,0) \subset (X,0)$ be a smooth subspace defined as the zero set of the analytic function
$g:(\C^{n+1},0) \to (\C^l,0)$, and let $\alpha \in S$ be as in remark \ref{ChoosingAlpha}-1. 
Then the pair $(X^0, V)$ satisfies Whitney's condition $a)$ at the origin if and only if the module 
$JM_g(f)_\alpha$ is contained in $JM(f)_\alpha^\dag$.
\end{Theorem}

\begin{corollary}\label{Restrictiontotherelativejacobianmodule}
    In the same setup of \ref{Whitagen}, let the smooth subspace $(V,0) \subset (X,0)$ be linear and
defined by the projection $g:(\C^{n+1},0) \to (\C^l,0)$ onto the first $l$ coordinates. If $h:(\C^{n+1},0)
\to (\C^{n+1-l},0)$ denotes the retraction over $(V,0)$, that is the projection onto the last $n+1-l$ coordinates,
then the pair $(X^0, V)$ satisfies Whitney's condition $a)$ at the origin if and only if the module $JM_g(f)_\alpha$
is contained in $JM_h(f)_\alpha^\dag$.
\end{corollary}
\begin{proof}
  Recall that 
   \[JM(f)_\alpha=\left< \left( \frac{\partial f_\alpha}{\partial z_0}\right), \cdots, 
                  \left( \frac{\partial f_\alpha}{\partial z_n}\right)\right>O_X^p\]
  where $\left(\frac{\partial f_\alpha}{\partial z_j}\right)= 
          \begin{pmatrix} \frac{\partial f_{\alpha_1}}{\partial z_j} \\ \vdots \\ \frac{\partial f_{\alpha_c}}{\partial z_j} 
               \end{pmatrix}$. Then, according to definition \ref{Reljacmoddef} we have that:
\[JM_g(f)_\alpha=\left< \left( \frac{\partial f_\alpha}{\partial z_l}\right), \cdots, 
                  \left( \frac{\partial f_\alpha}{\partial z_n}\right)\right>O_X^p\]
and
\[JM_h(f)_\alpha=\left< \left( \frac{\partial f_\alpha}{\partial z_0}\right), \cdots, 
                  \left( \frac{\partial f_\alpha}{\partial z_{l-1}}\right)\right>O_X^p\]
   Now, by definition, for a fixed map $(\phi,\psi):(\C,0) \to (X \times \check{\C}^p, 0)$, with $(\phi(\tau),
[\psi(\tau)]) \notin Z$ for $\tau \neq 0$, we have the ideal:
\begin{align*} I_\psi(JM(f)_\alpha \circ \phi)= & \left< \psi(\tau) \left( \frac{\partial f_\alpha}{\partial z_0}\circ \phi \right),
                          \cdots, \psi(\tau) \left( \frac{\partial f_\alpha}{\partial z_n}\circ \phi \right)
                                   \right>\C\{ \tau \} \\
                           = & \left< \tau^{r_0}w_0, \ldots, \tau^{r_{n}}w_{n}\right>\C\{\tau\}, \;
                         \mathrm{with}\; w_j \in \C\{\tau\} \; \mathrm{unit} \\
                           = & \left< \tau^k\right> \C\{\tau\}
\end{align*}
 But, by theorem \ref{Whitagen} we know that the pair $(X^0,V)$ satisfies Whitney's condition a) at 
the origin if and only if $JM_g(f)_\alpha \subset JM(f)_\alpha^\dag$. That is, for $j= l , \ldots, n$ we have that
\[\psi(\tau) \left( \frac{\partial f_\alpha}{\partial z_j}\circ \phi \right) \in \m_1 I_\psi(JM(f)_\alpha \circ \phi) =
      \left< \tau^{k+1} \right> \C \{\tau\}\]
so finally:
\begin{align*}
   \left< \tau^{r_0}w_0, \ldots, \tau^{r_{n}}w_{n}\right>\C\{\tau\} & = 
    \left< \tau^{r_0}w_0, \ldots, \tau^{r_{l-1}}w_{l-1}\right>\C\{\tau\} \\
   & = I_\psi(JM_h(f)_\alpha \circ \phi).
\end{align*}
 and the result follows.
\end{proof}

\begin{corollary}\label{AlphaWhitneyreljac} Let $\varphi:(\X,0) \to \C$  denote the specialization of $(X,0)$ to its
tangent cone $(C_{X,0},0)$, and let $\alpha \in S$ be as in remark \ref{ChoosingAlpha}-1. Then, the pair $(\X^0 ,Y)$ satisfies Whitney's condition $a)$ at the origin
if and only if $\frac{\partial F_\alpha}{\partial t} \in JM_\varphi (F)_\alpha^\dag$.
\end{corollary}
\begin{proof}
   For $(Y,0)\subset (\X,0) \subset (\C^{n+1} \times \C, 0)$ we have that the projection 
 \[\varphi:(\C^{n+1}\times \C,0) \to (\C,0)\]
 onto the last coordinate can be seen as the retraction over $(Y,0)$. Moreover, the subspace $(Y,0)$ is defined 
by the projection \[g:(\C^{n+1} \times \C,0) \to (\C^{n+1},0)\] onto the first $n+1$ coordinates, so the module
$JM_g(F)_\alpha=<\frac{\partial F_\alpha}{\partial t}> O_{\X}^p$, and the result follows from \ref{Restrictiontotherelativejacobianmodule}.
\end{proof}

\begin{remark}\label{FuncionaPuntosCerca}
    Proposition \ref{AlphaExplosionZ} gives us a relation between the blowup space $E_Z(\X \times \check{\P}^{c-1})$ and the limits
  of tangent hyperplanes for every point in a small enough neighborhood of the origin in $\X$. Since it is this
  relation what gives the key to derive  \ref{Limittangenthyperplane} to \ref{Restrictiontotherelativejacobianmodule}, these 
  results are also valid for every point in a small enough neighborhood of the origin in $\X$; all we have to 
  change is that the arcs $\phi: (\C,0) \to (\X^0, (z,t))$ arrive to the desired point. But more importantly, the 
  characterization of Whitney's condition a) given in corollary \ref{AlphaWhitneyreljac} is valid as stated for any 
  sufficiently close point $y \in Y$.  
\end{remark}

\section{The Main Theorem}

   Let $(X,0)$ be a reduced  germ of analytic singularity of pure dimension $d$, with reduced tangent cone
$C_{X,0}$, and let $(X,0)=\bigcup_{j=1}^r (X_j,0)$ be its irreducible decomposition. By lemma 
\ref{Especializaciones} $(\X,0)=\bigcup_{j=1}^r (\X_j,0)$ is the irreducible decomposition of the 
specialization space $\X$, where $(\X_j,0)$ is the specialization space of the irreducible component 
$(X_j,0)$ to its tangent cone $(C_{X_j,0},0)$. Moreover, if the germ $(X,0)$ doesn't have exceptional cones, 
 by lemma \ref{CompsIrrSinConosExcepcionales} the germs $(X_j,0)$ don't have exceptional cones
 either. These two results allow us to restrict ourselves to the case where \textbf{the germ $(X,0)$ 
 is irreducible} which we have been treating.\\
 
\begin{lema}\label{CompsIrrSinConosExcepcionales}
   The germ $(X,0)$ doesn't have exceptional cones if and only if for each 
$i \in \{1,\ldots,r\}$ the germ $(X_i,0)$ doesn't have exceptional cones. 
\end{lema}
\begin{proof}
   First of all, for a small enough representative of $X \subset \C^{n+1}$, we have the equality 
 $C(X)= \bigcup C(X_i)$ where $C(X_i)$ denotes the conormal space of the embedding $X_i \subset \C^{n+1}$, and 
 so the conormal map $\kappa_{X_i}$ is equal to the restriction of $\kappa_X$ to $C(X_i)$.
 Moreover, we know that the strict transform $ \overline{e_0^{-1}( X_i \setminus \{0\})}$ is equal to 
 the blowing-up $E_0X_i \to X_i$, and since for every arc $\phi:(\C,0) \to (X,0)$ there exists a $j \in \{1,\ldots,r\}$
 such that $\phi$ factorizes through $X_j$, we have the equality  $\P C_{X,0}= \bigcup \P C_{X_i,0}$. All these
 imply that for each $i \in \{1,\ldots,r\}$ the normal conormal diagram
  \[\xymatrix{E_0C(X_i)\ar[r]^{\hat{e}_0}\ar[dd]^{\kappa'_{X_i}}\ar[ddr]^\zeta & C(X_i)\ar[dd]^{\kappa_{X_i}} \\
             &  \\
             E_0X_i\ar[r]_{e_0}  & X_i }\]
  is canonically embedded in the normal conormal diagram of $X$:
 \[\xymatrix{E_0C(X)\ar[r]^{\hat{e}_0}\ar[dd]^{\kappa'_X}\ar[ddr]^\zeta & C(X)\ar[dd]^{\kappa_X} \\
             &  \\
             E_0X\ar[r]_{e_0}  & X }\]
  
   Now, the germ $(X,0)$ doesn't have exceptional cones if and only if every irreducible 
 component $W_\alpha$ of the fiber $|\kappa_X^{-1}(0)|= \bigcup |\kappa_{X_i}^{-1}(0)|$ is equal to the 
 projective dual of an irreducible component $V_\alpha$ of the tangent cone $\P C_{X,0}$, that is an 
 irreducible component of one of the tangent cones $\P C_{X_i,0}$. Finally, since for a reduced 
 projective subvariety the double dual $\check{\check{Y}}$ is equal to $Y$, then two projective 
 subvarieties $Y_1$ and $Y_2$ of $\P^n$ are different if and only if their duals are different $\check{Y_1} \neq
 \check{Y_2}$. This prevents the appearance of a possible exceptional cone of $X_j$ having the same dual
 as an irreducible component of $\P C_{X,0}$ which finishes the proof.
\end{proof} 

    As we have said before, the first step in our objective of constructing a Whitney stratification of 
 $(\X,0)$ having the parameter axis $(Y,0)$ as a stratum, is proving that the pair $(\X^\circ,Y)$ 
 satisfies Whitney conditions $a)$ and $b)$ at the origin.  Since we are assuming $(\X,0)$ irreducible, what
 we have to prove, according to corollary  \ref{AlphaWhitneyreljac} is that for an
 $\boldsymbol\alpha$ \textbf{chosen as in remark \ref{ChoosingAlpha}-1, which we will fix from this point on}, 
 $\frac{\partial F_\alpha}{\partial t} \in JM_\varphi(F)_\alpha^\dag$.
 So in terms of \ref{Moduloseidealjacobiano}, what we must prove (assuming we know that the rank
 of the $\alpha-$relative jacobian module is the codimension $c$) is that every minor $M$ in 
 $J_c(JM_\alpha(F))$ depending on $\frac{\partial F_\alpha}{\partial t}$ satisfies 
 $M \in J_c(JM_\varphi(F)_\alpha)^\dagger$.   
 We will prove this using \ref{DependenciaMeromorfa}, and since we are working with ideals, it
 leads us to consider the normalized blowup of $\X$ along the ideal $J_c(JM_\varphi (F)_\alpha)$. 
 Moreover, by remark \ref{AlphaNash}, the blowup of $\X$ along the ideal $J_c(JM_\varphi (F)_\alpha)$
 gives the relative Nash modification $\nu_\varphi:\Na_\varphi \X \to \X$.\\

  \begin{lema}\label{rankofJM(F)} 
   The $\alpha$-Jacobian module $JM(F)_\alpha$ has rank $c$ on $(\X,0)$. 
 \end{lema}
  \begin{proof}
     By definition, the rank of a module over the integral domain
  $O_{\X,0}$ is the dimension as a vector space over the quotient field $Q(O_{\X,0})$ of the
  vector space $Q(O_{\X,0}) \bigotimes JM(F)_\alpha O_{\X,0}$. 

    Consider the presentation
 \[ O_{\X,0}^r \longrightarrow O_{\X,0}^{n+2} \stackrel{[DF]_\alpha}{\longrightarrow} JM(F)_\alpha O_{\X,0} \to 0\]
   where $[DF]_\alpha$ denotes the jacobian matrix of the map $F_\alpha:\C^{n+2} \to \C^c$, which defines this map.
  By tensorizing this sequence by the field $Q(O_{\X,0})$, we obtain the sequence
 \[ Q(O_{\X,0})^r \longrightarrow Q(O_{\X,0})^{n+2} \stackrel{[DF]_\alpha}{\longrightarrow} 
   Q(O_{\X,0}) \bigotimes JM_\alpha(F)O_{\X,0} \to 0\]
  where the map defined by the jacobian matrix remains surjective. Remark that we now have
  that the rank of the module $JM_\alpha(F)O_{\X,0}$ is equal to the rank of the matrix $[DF]_\alpha$ when
  considering its entries as members of the quotient field $Q(O_{\X,0})$.\\

    Our choice of $\alpha$ garanties the existance of a non zero $c \times c$ minor in $O_{\X,0}$. This implies
  that the ideal $J_c(JM(F)_\alpha O_{\X,0})$ of $O_{\X,0}$ generated by all the $c \times c$ minors of 
  the matrix $[DF]_\alpha$ is different from zero. Moreover
  since the matrix $[DF]_\alpha$ is of size $c \times (n+2)$, then the ideal $J_{c+1}(JM(F)_\alpha O_{\X,0})$ is 
  equal to the zero ideal. This remains true when considering the minors as elements of the quotient 
  field $Q(O_{\X,0})$, and so the rank of the matrix $[DF]_\alpha$ is equal to $c$ which finishes the proof.
  \end{proof}

     We know that the pair $(\X^0,Y)$ satisfies Whitney's conditions at every point $y$ of $Y$ with the 
  possible exception of the origin, so we have by \ref{FuncionaPuntosCerca} that every minor $M$ in 
  $J_c(JM(F)_\alpha)$ depending on $\frac{\partial F_\alpha}{\partial t}$ satisfies 
  $M \in J_c(JM_\varphi (F)_\alpha)^\dagger$ in $O_{\X,y}$ for all these points. What we are going to prove 
  in proposition \ref{PruebamayorqueunoCI} is that this condition carries over to the origin under  
  the assumption that $(X,0)$ does not have exceptional cones.
  
\begin{remark}\label{EspacioZIntCom} $\,$
  \begin{enumerate}
\item  The fact proven in proposition \ref{IsomorfismosConormales}, that the isomorphism between the conormal 
       space $C(\X\setminus \X(0))$ and $C(X) \times \C^*$ is given by a natural projection implies
       that the vertical hyperplane $\{t=0\}:=[0:\cdots:0:1] \in \check{\P}^{n+1}$ is not 
       tangent to any point $(z,t) \in \X \setminus \X(0)$. This is equivalent, by corollary 
       \ref{NoHayHiperplanoVertical}, to $\frac{\partial F_\alpha}{\partial t} \in 
        \overline{JM_\varphi (F)_\alpha}$ in $O_{\X,(z,t)}$ for every point $(z,t) \in \X \setminus \X(0)$.
  
\item  When $(\X,0)$ is a complete intersection, the center of the blowup defined by the ideal 
     $J_c(JM_\varphi (F))$ is set theoretically the relative singular locus of $\X$. Moreover, since
     in this case, the tangent cone $(C_{X,0},0)$ is a complete intersection,   
     the equality $\frac{\partial F_i}{\partial z_j}(z,0) = \frac{\partial f_{m_i}}{\partial z_j}(z)$ 
     give us that the restriction of the ideal $J_c(JM_\varphi (F))$ to the special fiber is equal to the 
     jacobian ideal $J_{C_{X,0}}$ of the tangent cone $C_{X,0}$ in $O_{C_{X,0}}$. This implies that the strict 
     transform of $\X(0)$ with respect to this blowup is equal to the Nash modification $\Na C_{X,0}$ of the fiber.

\item  Even though we are considering that $(X,0)$ and as a result $(\X,0)$ are \linebreak irreducible germs, this
  doesn't mean that the tangent cone $(C_{X,0})$ is irreducible. The problem with this is that the restriction 
  of the ideal $J_c(JM_\alpha(F))$ to the special fiber $\X(0)$ may vanish in an irreducible component
  of the tangent cone $(C_{X,0},0)$ and so its strict transform will no longer be equal to the Nash 
  modification $\Na C_{X,0}$. 
  \end{enumerate}
\end{remark}

\begin{lema}\label{BuenidealI}
     For a reduced and irreducible germ $(X,0)$ of analytic singularity with reduced tangent cone $(C_{X,0},0)$, 
  there exists an ideal $I \subset O_{\X,0}$  such that:
  \begin{enumerate}
   \item The analytic subset $V(I) \subset \X$ defined by $I$ contains the relative singular locus 
         $\mathrm{Sing}_\varphi \X := \bigcup_t  \mathrm{Sing} \X(t)$.
   \item The blowup of $\X$ along $I$ is equal to the relative Nash modification of $\X$, that is
         $E_I\X \cong \Na_\varphi \X$. 
   \item  The blowup of the special fiber $\X(0)$ along the ideal 
         $IO_{\X(0),0}$ defined by the restriction of $I$ to $\X(0)$ is isomorphic to the Nash
         modification $\Na C_{X,0}$. 
  \end{enumerate} 
\end{lema}
\begin{proof}
      Let $F: (\C^{n+1} \times \C,0) \to (\C^p,0)$ denote the germ of analytic map defined by the $p$ series 
 $F_1, \ldots, F_p \in \C\{z_0,\ldots,z_n,t\}$, such that $(\X,0)=(F^{-1}(0),0)$. Let $[D_\varphi F]$ denote
 the relative jacobian matrix, and define the $p \times (n+1)$ matrix $A$ by setting the $t$ coordinate
 to $0$, that is $A=[D_\varphi F (z,0)]$. By definition, $A$ is the jacobian matrix of the map
 $g: (\C^{n+1},0) \to (\C^p,0)$ defined by the homogeneous polynomials $g_i=F_i(z,0)$ such that
 $(C_{X,0},0)=(g^{-1}(0),0)$. Let $c$ be the codimension of $\X$ in $\C^{n+1} \times \C$, then $c$ is also
 the codimension of $C_{X,0}$ in $\C^{n+1}$, and let $S$(resp. $S'$) denote the set of increasing sequences of 
 $c$-positive integers less than $p+1$ (resp. $n+2$). For $\alpha=(\alpha_1, \ldots, \alpha_c) \in S$, 
 and $\beta=(\beta_1, \ldots, \beta_c)\in S'$, $g^{\alpha \beta}$ will denote the minor of $A$ obtained 
 by considering the rows determined by $\alpha$ and the columns determined by $\beta$.\\

    Let $C_{X,0}=\bigcup_{j=1}^l V_j$ be the irreducible decomposition of the tangent cone. By
  the proof of \ref{NashNobileCoro} there exist
  $\alpha^1, \ldots, \alpha^l$ in $S$ and functions $h_1, \ldots, h_l \in O_{C_{X,0},0}$, with 
  $h_i=0$ on $\bigcup_{j \neq i} V_j$ and $h_i \neq 0$ on $V_i$, such that the blowup of $C_{X,0}$ along
  the ideal $J=\left<\sigma_\beta:= \sum_{i=1}^l h_ig^{\alpha^i,\beta}, \;\beta \in S'\right>$ gives 
  the Nash modification $\Na C_{X,0}$.\\

    Now, since for each $\alpha^i$ there is a non-zero minor of the matrix $[Dg]_{\alpha^i}$, 
  the corresponding minor of the matrix $[D_\varphi F]_{\alpha^i}$ is not identically zero. Since 
  by hypothesis $\X$ is irreducible then the proof of \ref{NashNobileCoro} tells us that this condition
  is enough to prove that the blowup of $\X$ along the ideal $J_c(JM_\varphi(F)_{\alpha^i})$ gives
  the relative Nash modification $\Na_\varphi \X$.\\

   Let $F^{\alpha \beta}$ denote the minor of $[D_\varphi F]$ obtained by considering the rows determined by 
   $\alpha$ and the columns determined by $\beta$, and define the ideal $I=\left<\rho_\beta:= 
   \sum_{i=1}^l h_i F^{\alpha^i,\beta}, \;\beta \in S'\right>$ where the $h_i$'s are the same we used
   for the tangent cone. Now, by construction, the blowup of the special fiber $\X(0)$ along the ideal 
   $IO_{\X(0),0}$ is isomorphic to the Nash modification $\Na C_{X,0}$, and since for any point $(z,t)$
   in the relative singular locus all the $c \times c$ minors of $[D_\varphi F]$ vanish, then 
   we have the inclusion $\mathrm{Sing}_\varphi \X \subset V(I)$. All that is left to prove, is that
   the blowup of $I$ gives $\Na_\varphi \X$. \\

      Let $x= (z,t)$ be a point in the relative smooth locus of $\X$ and $T_x\X(t)^0=[a_0: \ldots :a_N]$ denote
  the coordinates of the point of $\P^N$ corresponding to the direction of the tangent space to the fiber $\X(t)$ 
  at $x$ by the Plucker embedding of the grassmannian $G(d,n+1)$ in the projective space $\P^N$. If $(z,t)$ is
  sufficiently general then for each of the $\alpha^i$'s we have:
  \[[F^{\alpha^i,\beta^0}: \cdots : F^{\alpha^i,\beta^N}]=[a_0: \ldots :a_N]\]
  where we have ordered the $\beta$'s lexicographically. This means that there exist $\lambda_1, \ldots, \lambda_l
  \in \C$ such that for every $\alpha^1, \ldots, \alpha^l $ and  $\beta^k \in S'$ we have:
  \[ F^{\alpha^i,\beta^k}= \lambda_i a_k\]
  which implies that for each $\beta^k \in S'$:
  \[\rho_{\beta^k}(x)= \sum_{i=1}^l h_i F^{\alpha^i,\beta^k}(x) = \sum_{i=1}^l h_i \lambda_i a_k =
     a_k\sum_{i=1}^l \lambda_i h_i \]
  and so $[\rho_\beta (x)]= [a]$ in $\P^N$. Finally, since the $\lambda's$ are non zero constants,
  the function $\sum_{i=1}^l \lambda_i h_i$ can not be identically zero. This implies that the 
  equation $[\rho_\beta (x)]= [a]$ in $\P^N$ is true for every point $x$ in an open dense set $U \subset \X$
  which finishes the proof.   
\end{proof}

\begin{proposition}\label{VerticalComponentsCI} 
   Let $\nu_\varphi: \Na_\varphi \X \to \X$ be the relative Nash modification of $\varphi: \X \to \C$.
Let $Z \subset \X$ be the subspace defined by the ideal $I$ of \ref{BuenidealI}, and let $D$ be
the divisor defined by $I$ in $\Na_\varphi \X$, that is $D=\nu_\varphi^{-1}(Z)$. If the germ $(X,0)$ does not
have exceptional cones, then $\nu_\varphi^{-1}(Z\setminus \X(0))$ is dense in $D$. That is, the exceptional 
divisor $D$ of $\Na_\varphi \X$ does not have \textbf{\emph{vertical components}} over $\X(0)$.
\end{proposition}
\begin{proof}
    We know that $\X(0)$ is isomorphic to the tangent cone 
   $C_{X,0}$. Now, by \ref{BuenidealI} the strict transform of $\X(0)$ in $\Na_\varphi(\X)$ is isomorphic 
   to the Nash modification $\upsilon:\Na C_{X,0} \to C_{X,0}$. Moreover, by the definition of blowup, 
   $\upsilon^{-1}(Z\cap \X(0))$ is a divisor (of dimension $d-1$).\\
    
    Now, by \ref{NashrelativevsNash} if $(z,0,T) \in \nu_\varphi^{-1}(z,0) \subset \Na_\varphi \X$ then 
    the $d-$plane $T$ is via $\Gamma$ a limit of tangent spaces to $X$ at $0$, that is the point 
    $(0,T) \in \nu^{-1}(0) \subset \Na X$. But, since by hypothesis the germ $(X,0)$ does not have
    exceptional cones , then $T$ is tangent to the tangent cone $C_{X,0}$.\\ 
    
       We want to prove that the total transform $\nu_\varphi^{-1}(\X(0))$ coincides with 
    the strict transform $\Na C_{X,0}$, that is, we need to prove that the point $(z,0,T)$ is in 
    $\Na C_{X,0}$. For this purpose all that is now left to prove is that $T$ is tangent to $C_{X,0}$ 
    at the point $p=(z)$. \\

        Let $\delta:(\C, \C \setminus \{0\},0) \to (X, X^0, 0)$ be an arc such that its lift  
   $\widetilde{\delta}$ to $\Na X$ has the point $(0,T) \in \Na X$ as endpoint.
   
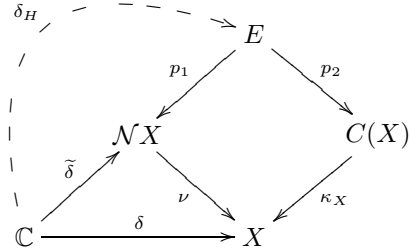
\[\xymatrix{
        &  &  E \ar[dl]_{p_1} \ar[dr]^{p_2}& \\
       & \Na X \ar[dr]_{\nu} & & C(X) \ar[dl]^{\kappa_X}\\
          \C \ar[rr]^{\delta} \ar[ur]^{\widetilde{\delta}} \ar @/^4.5pc/ @{-->}[uurr]^{\delta_H}& & X &
     }\]
    By construction $\delta(\C\setminus\{0\})$ is contained in the smooth locus $X^\circ$, and if
    we denote by $E^\circ$ the inverse image $p_1^{-1}(\nu^{-1}(X^0))$ , then by \ref{ConormalvsNash} 
    the open subset $E^\circ$ is dense in $E$, and it defines a locally trivial fiber bundle over 
    $\nu^{-1}(X^\circ)$. This implies that for any point $(0,T,H) \in p_1^{-1}(0,T)$ the arc 
    $\widetilde{\delta}$ can be lifted to an arc $\delta_H$ 
    having the point $(0,T,H)$ as endpoint. So now we have transformed the problem into proving
    that any hyperplane $H \in \check{\P}^n$, such that $T \subset H$, is a tangent hyperplane to
    $C_{X,0}$ at the point $p=z$. \\
    
     Going back again to the diagram of \ref{NashrelativevsNash}:
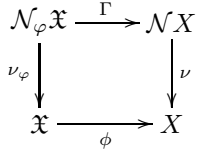
    \[\xymatrix{ \Na_\varphi \X \ar[r]^\Gamma \ar[d]_{\nu_\varphi} &  \Na X \ar[d]^{\nu}\\
                   \X  \ar[r]_\phi & X }\]   
     we have that for any sequence $\{(z_m,t_m)\}$ in the smooth part of $\X \setminus 
    \X(0)$ tending to the point $(z,0)$ in the special fiber $\X(0)$, we have a corresponding 
    sequence $\{(t_m z_m)\}$ tending to the origin in $X$. The final step of the proof
    is now a consequence of the projective duality obtained from the normal/conormal diagram:
     
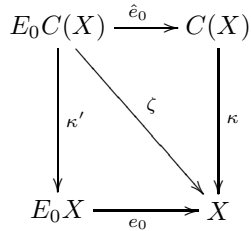
\[\xymatrix{E_0C(X)\ar[r]^{\hat{e}_0}\ar[dd]^{\kappa'}\ar[ddr]^\zeta & C(X)\ar[dd]^\kappa
            \\
             & &  \\
             E_0X\ar[r]_{e_0}  & X }\]
     since the sequence $\{t_m z_m\} \subset X \setminus \{0\}$ tending to the origin gives us the
     sequence $\{(t_m z_m), [z_m]\}$ in $E_0X$, the blowup of $X$ at $0$, which tends to the 
     point $(0,[z])$ in the exceptional divisor $\P C_{X,0}$. In the same way, we obtain the sequence
     $\{(t_m z_m, [z_m], H_m)\}$ in $E_0C(X) \subset X \times \P^n \times \check{\P}^n$ tending 
     to the point $(0, [z],H)$ in $G = \zeta^{-1}(0)$. Recall that
     if $|G|= \bigcup_\alpha G_\alpha$ is the irreducible decomposition of the reduced space
     $|G|$, then each $G_\alpha$ is the conormal space of an irreducible component of $\P C_{X,0}$.
     To finish the proof, note that so far we have proved that $\nu_\varphi^{-1}(\X(0))$ 
     is just $\Na C_{X,0}$ and so $\nu_\varphi^{-1}(Z(0))$ is of dimension $d-1$, whereas an 
     irreducible component of $D$ is of dimension $d$.
\end{proof}

\begin{corollary}\label{DimensionCorrecta}
   Let $\mathrm{Sing} \X(0)$ denote the singular locus of the special fiber, then the dimension of 
  $\nu_\varphi^{-1}(\mathrm{Sing} \X(0))$ is less or equal than $d-1$.
\end{corollary}
\begin{proof}
    By definition of the ideal $I$, the analytic subset $\mathrm{Sing} \X(0)$ is contained in the 
 subspace $Z$ defined by $I$. Then we have the inclusion $\nu_\varphi^{-1}(\mathrm{Sing} \X(0)) 
 \subset \nu_\varphi^{-1}(Z(0))$ and 
 by proposition \ref{VerticalComponentsCI} the dimension of $\nu_\varphi^{-1}(Z(0))$ is equal to $d-1$
 which finishes the proof.
\end{proof}

   Note that the following result does not uses the irreducible hypothesis, and so is valid
  in a more general setting.

\begin{lema}\label{LemaImagenInversaProducto}
   Let $Y$ denote the smooth subspace $0 \times \C \subset \X$ as before, let $\nu: \Na X \to X$ be
 the Nash modification of $X$, and let $\widetilde{\nu_\varphi}:\widetilde{\Na_\varphi \X} \to \X$ be the normalized relative 
 Nash modification of $\X$. Then:
 \begin{enumerate}
  \item  If the germ $(X,0)$ doesn't have exceptional cones we have the set-theoretical equality:
 \[\left |\nu_\varphi^{-1}(Y)\right|= \left|Y \times \nu^{-1}(0)\right|\]
  \item The set theoretical inverse image $|\widetilde{\nu_\varphi}^{-1}(Y\setminus \{0\})|$ is dense in $|\widetilde{\nu_\varphi}^{-1}(Y)|$. 
 \end{enumerate}
\end{lema}
\begin{proof}
   From proposition \ref{NashrelativevsNash} we have the commutative diagram:
  \[\xymatrix{ \Na_\varphi \X \ar[r]^\Gamma \ar[d]_{\nu_\varphi} &  \Na X \ar[d]^{\nu}\\
                   \X  \ar[r]_\phi & X }\]   
  where $\phi$ and $\Gamma$ are surjective. The morphism $\phi$ is the restriction
  to $\X$ of the map $\C^{n+1} \times \C \to \C^{n+1}$ defined by $(z_0, \ldots,z_n,t) \mapsto
  (tz_0, \ldots, tz_n)$ which is an isomorphism on $\C^{n+1} \times \C^*$. This implies in particular
  that the restriction of the differential $D\phi$ to the tangent space $T_{(z,t)}\X(t)$ maps it 
  isomorphically to $T_{(tz)}X$, where $(z,t)$ is a smooth point of the fiber $\X(t)$ with $t \neq 0$.
  But the restriction of $D\phi$ to $T_{(z,t)}\X(t)$ is $t$ times the identity $\mathrm{Id}$, which implies
  that $\nu_\varphi^{-1}(Y\setminus \{(0,0)\}= Y\setminus \{(0,0)\} \times \nu^{-1}(0))$ and as a consequence
  $\nu_\varphi^{-1}(0,0)$ contains $\nu^{-1}(0)$. Finally, from the proof of proposition \ref{VerticalComponentsCI}
  we know that the fiber $\nu_\varphi^{-1}(\X(0))$ is equal to the Nash modification of the tangent 
  cone $C_{X,0}$, so the fiber $\nu_\varphi^{-1}(0,0)$ is equal to the set of limits of tangent spaces to 
  $C_{X,0}$ which coincides with $\nu^{-1}(0)$ since the germ $(X,0)$ doesn't have exceptional cones.\\

    To prove $2)$, note that since $\nu_\varphi^{-1}(Y)$ has a product structure we already have that 
  $\nu_\varphi^{-1}(Y\setminus \{0\})$ is dense in $Y$, and so we need to study how the normalisation
  $n:\widetilde{\Na_\varphi \X} \to \Na_\varphi \X$ affects this subspace. Let $(0,0,T) \in \Na_\varphi \X$
  be a point over the origin in $\X$. Since by hypothesis $\X$ is irreducible, the space $\Na_\varphi \X$ is 
  also irreducible, however it may not be locally irreducible so the germ $(\Na_\varphi \X,(0,0,T))$ may have
  an irreducible decomposition of the form $(\Na_\varphi \X,(0,0,T))= \bigcup_j (W_j,(0,0,T))$. Now, by 
  \cite[Section 4.4]{DeJo-Pfi}, we have that the normalisation map is finite, and over 
  $(\Na_\varphi \X,(0,0,T))$ in the normalised space 
  $\widetilde{\Na_\varphi \X}$ we have a multigerm $\bigsqcup_j (\widetilde{W_j},p_j) $ such that:
  \begin{enumerate}
   \item The germ $(\widetilde{W_j},p_j)$ is irreducible, and corresponds to the normalisation of  
         $(W_j,(0,0,T))$.
   \item For every $j$ we have that $n^{-1}(0,0,T) \cap \widetilde{W_j}=\{p_j\}$. 
  \end{enumerate}
    This implies that if $\nu_\varphi(\nu_\varphi^{-1}(Y) \cap W_j)=Y$, then set-theoretically $\widetilde{\nu_\varphi}^{-1}(Y\setminus 
  \{0\}) \cap \widetilde{W_j}$ is dense in $\widetilde{\nu_\varphi}^{-1}(Y) \cap \widetilde{W_j}$, and so all we have to prove
   is that every $W_j$ satisfies this condition.\\

     Since the open set of relative smooth points $\X_\varphi^0 \setminus \X(0)$ is dense in $\X$, then its preimage 
  $\nu_\varphi^{-1}(\X_\varphi^0\setminus \X(0))$ is dense in $\Na_\varphi \X$ and so it intersects every irreducible 
  component $W_j$ in an open dense set $U_j$. This means that there exists an arc contained in $U_j$
   \begin{align*}
           \mu:(\C, \C\setminus \{0\},0) &\to (W_j, U_j, (0,0,T))\\
                     \tau & \mapsto \left( z(\tau), t(\tau), T(\tau)\right)   
   \end{align*}

   having $(0,0,T)$ as endpoint; moreover by composing it with $\nu_\varphi$ we get an arc 
   \[\widetilde{\mu}:(\C, \C \setminus \{0\},0) \to (\X,\X_\varphi^0 \setminus \X(0),(0,0))\] 
   contained in $\X_\varphi^0 \setminus \X(0)$ having the origin as endpoint. \\

     Let $\widetilde{\mu}=(z(\tau),t(\tau))$ and let $\alpha \in \C^*$, 
   by propositions \ref{IsomorfismosConormales} and \ref{NashrelativevsNash}, this arc can be "verticalized"
   to an arc $\widetilde{\mu_\alpha}:(\C, \C\setminus\{0\},0) \to \left(\X(\alpha), \X^0(\alpha), (0, \alpha)\right)$ 
   as follows:
   \begin{align*}
    (\C, \C\setminus\{0\},0) &\to 
     \left( \X,\X_\varphi^0 \setminus \X(0),(0,0) \right) \longrightarrow (X,X^0,0) \longrightarrow 
    \left(\X(\alpha), \X(\alpha)^0, (0, \alpha)\right) \\ \tau &\mapsto
    \left(z(\tau), t(\tau)\right) \longmapsto \left( t(\tau)z(\tau) \right) \longmapsto
    \left( \frac{ t(\tau)z(\tau)}{\alpha}, \alpha \right)
\end{align*}
     Since the canonical isomorphism between two fibers $\X(\alpha_1)$ and $\X(\alpha_2)$ used
   here is given by $(z,\alpha_1) \mapsto (\frac{\alpha_1 }{\alpha_2}z,\alpha_2)$, for 
   every smooth point the tangent map acts as $\frac{\alpha_1}{\alpha_2}$ times the identity on 
   the embedded tangent space leaving it invariant. Now, since the arc is contained in the
   smooth locus $\X^0(\alpha)$ it has a unique lift to an arc 
 \[\mu_\alpha:(\C,\C\setminus\{0\},0) \to (\Na_\varphi \X, \nu_\varphi^{-1}(\X_\varphi^0),(0,\alpha,T))\]
  having as endpoint the point $(0,\alpha,T)$. Moreover for every $\tau_0$ close enough to the origin in 
  $\C$ the point $(z(\tau_0), t(\tau_0), T(\tau_0))$ is in $W_j$ and since the arc $\mu_{t(\tau_0)}$ 
  passes through this point, then it is completely contained in $W_j$, in particular the 
  endpoint $(0,t(\tau_0),T)$ is in $W_j$ which finishes the proof.     
\end{proof}

   We are now in position to prove that $\frac{\partial F_\alpha}{\partial t}$ is strictly dependent on 
 $JM_\varphi (F)_\alpha$ at $0$.
 
 \begin{proposition}\label{PruebamayorqueunoCI}
     If the germ $(X,0)$ does not have exceptional cones then every minor $M$ in 
 $J_c(JM(F)_\alpha)$ depending on $\frac{\partial F_\alpha}{\partial t}$ satisfies $M \in 
 J_c(JM_\varphi (F)_\alpha)^\dagger$ in $O_{X,0}$.   
  \end{proposition}
  \begin{proof}
      Let $M$ be a minor in $J_c(JM(F))$ that depends on $\frac{\partial F_\alpha}{\partial t}$, and let 
    $W \subset \X$ be the subspace defined by the ideal $J_c(JM_\varphi (F)_\alpha)$. 
    Note that by definition, not only the $t$-axis $Y$, but the entire relative singular locus 
    $\mathrm{Sing}_\varphi \X$ is contained in $W$. Let $\widetilde{\nu_\varphi}:\widetilde{\Na_\varphi \X} \to \X$ be the 
    normalized blowup of $\X$ along $J_c(JM_\varphi (F)_\alpha)$, and let $\overline{D}$ be its exceptional divisor. 
    By considering a small enough neighborhood of the origin in $\X$, or in other words a small enough representative
    of the germ $(\X,0)$ we can assume that the divisor $\overline{D}$ has a finite number of  irreducible 
    components, and every irreducible component of $\overline{D}$ intersects $\widetilde{\nu_\varphi}^{-1}(0)$.
    Thanks to the fact that each irreducible component $\overline{D}_k$ is mapped by the normalisation map
     $n: \widetilde{\Na_\varphi \X} \to \Na_\varphi(\X)$ to an irreducible component $D_j$ of $ D = 
     |\nu_\varphi^{-1}(W)|$ these conditions are also verified in $\Na_\varphi(\X)$. \\

       Let $b \in \overline{D}$ be a point in the exceptional divisor lying over $W(0)$. Now, since 
    $\overline{D}$ is a divisor, the ideal $J_c(JM_\varphi (F)_\alpha) \circ \widetilde{\nu_\varphi}$ is locally invertible, so at each 
    $b \in \overline{D}(0)$ it is generated by a single element $g \circ \widetilde{\nu_\varphi}$, where $g \in J_c(JM_\varphi (F)_\alpha)$. 
    By proposition \ref{DependenciaMeromorfa}, we need to prove that for every such $b$ the function 
    $M \circ \widetilde{\nu_\varphi}$ lies in the product $I(Y,\overline{D}_k)J_c(JM_\varphi (F)_\alpha) \circ \widetilde{\nu_\varphi}$, or equivalently
    (from the proof of the proposition) that the meromorphic function $k$ locally defined by $\frac{M \circ \widetilde{\nu_\varphi}}
    {g \circ \widetilde{\nu_\varphi}}$ is holomorphic and vanishes at $b$ if $b$ lies over $(0,0) \in Y$.\\ 

       Note that if $\widetilde{\nu_\varphi}(b)$ is not in $Y$ then the ideal $I(Y,\overline{D}_k)O_{\widetilde{\Na_\varphi \X},b}$
    is not a proper ideal and so all we need to prove is that $M \circ \widetilde{\nu_\varphi}$ belongs to the ideal $J_c(JM_\varphi 
    (F)_\alpha) \circ \widetilde{\nu_\varphi}$, which by proposition \ref{DependenciaMeromorfa} is equivalent to $k$ being holomorphic
    and also to $M \in \overline{J_c(JM_\varphi (F)_\alpha)}$. Now, by remark \ref{EspacioZIntCom}-1, for any 
    point $(z,t) \in \X\setminus \X(0)$ we already have $M \in \overline{J_c(JM_\varphi (F)_\alpha)}$ which implies
    that the function $k$ is holomorphic on $\overline{D} \setminus \overline{D}(0)$, and so its polar locus is 
    contained in $\overline{D}(0)$.\\

       Let $(z,0) \in W$ such that $(z,0)$ is not in $\mathrm{Sing}_\varphi \X$, that is $(z,0)$ is a smooth
    point of both the space $\X$ and the special fiber $\X(0)$. Then, the vertical hyperplane 
   $H=[0:\cdots:0:1] \in \check{\P}^{n+1}$ cannot be tangent to $\X$ at $(z,0)$ and so by remark \ref{EspacioZIntCom}-1 
   we have $M \in \overline{J_c(JM_\varphi (F)_\alpha)}$ and $k$ holomorphic. Indeed, if $H$ is tangent to 
   $\X$ at the point $(z,0)$, then the point $(z,0)$ is a singular point of $\X \cap H= \X(0)$. This implies 
   that the polar locus of $k$ is contained in $\widetilde{\nu_\varphi}^{-1}(\mathrm{Sing}\X(0))$, but by corollary 
   \ref{DimensionCorrecta} the dimension of $\nu_\varphi^{-1}(\mathrm{Sing}\X(0))$ is less than or equal
   to $d-1$, and since the normalisation map is finite we also have dim $\widetilde{\nu_\varphi}^{-1}(\mathrm{Sing}\X(0))< d$,
   that it has codimension at least $2$.
   However, in a normal space the polar locus of a meromorphic function is of codimension 1 or empty 
   (\cite[Thm. 71.12, p. 307]{Ka}), which implies that $k$ is holomorphic at every point $b \in \overline{D}$.\\
   
       All that is left to prove is that the holomorphic function $k$ vanishes at every point $b \in \overline{D}$ 
  lying over $Y$. Since for any point $y \neq 0 \in Y$ the pair $(\X^0,Y)$ satisfies Whitney's condition a) 
  at $y$ we have that $k$ vanishes on $\widetilde{\nu_\varphi}^{-1}(Y\setminus \{(0,0)\})$, and by continuity it vanishes on its 
  closure in $\widetilde{\Na_\varphi \X}$. But by lemma \ref{LemaImagenInversaProducto}-2 the aforementioned 
  closure is equal to $\widetilde{\nu_\varphi}^{-1}(Y)$, and so we have that the function $k$ vanishes at 
  any point $b$ lying over $(0,0) \in Y$.  
  \end{proof} 
 
    Let $Z \subset \X$ be the subspace defined by the ideal $I$ of \ref{BuenidealI} as before.
  Note that the key point in proving the previous proposition is the inequality
  dim $\nu_\varphi^{-1}(\mathrm{Sing} \X(0)) < d$ which was a consequence of \ref{VerticalComponentsCI} 
  and this gives us the following result. 
    
  \begin{proposition}\label{EquivSinConosSinComps}
     Let $(X,0)\subset (\C^{n+1},0)$ be a reduced and irreducible $d$ dimensional germ of analytic singularity
    such that the tangent cone is reduced. Then 
     $(X,0)$ does not have exceptional cones if and only if $\nu_\varphi^{-1}(Z)$ does not have vertical 
     components over $\X(0)$.
  \end{proposition}
  \begin{proof}
     If $(X,0)$ does not have exceptional cones, then it is proposition \ref{VerticalComponentsCI}.
     On the other hand, if $\nu_\varphi^{-1}(Z)$ does not have vertical components over $\X(0)$
     then corollary \ref{DimensionCorrecta} and the proof of proposition \ref{PruebamayorqueunoCI} gives us that 
     the pair $(\X^0,Y)_0$ satisfies Whitney's condition a) at the origin, and by \ref{PropoGeneral*} this 
     is equivalent to $(\X,0)$ having no exceptional cones. Finally, this implies that $(X,0)$ does not have
     exceptional cones either.
  \end{proof}

\begin{remark}\label{Remarkexcepcional}
     Note that if $(X,0)$ has exceptional cones then, $(\X,0)$ also has exceptional cones.
  \end{remark}

     Indeed, if $\kappa_\X:C(\X) \to \X$ is the conormal space of $\X$ and $\kappa_X: C(X) \to X $ the 
   conormal space of $X$, then $\kappa_\X^{-1}(Y\setminus \{0\}) = Y \setminus \{0\} \times \kappa_X^{-1}(0)$
   and so $ \kappa_\X^{-1}(Y)$ contains $Y \times \kappa_X^{-1}(0)$. In particular, if $H=[a_0: \cdots 
   a_n] \in \kappa_X^{-1}(0) \subset \check{\P}^n$, but $H$ is not tangent to the tangent cone $C_{X,0}$,
   then $\widetilde{H}=[a_0:\cdots:a_n:0] \in \kappa_\X^{-1}(0)\subset \check{\P}^{n+1}$ and it can 
   not be tangent to the tangent cone $C_{\X,0}= C_{X,0} \times \C$. \\

We can summarize all we have done so far with the following theorem:
 
 \begin{Theorem}\label{Equivalenciashipersuperficies}
   Let $(X,0)$ be a reduced and equidimensional germ of complex analytic singularity, and suppose
   that its tangent cone $C_{X,0}$ is reduced. Then the following 
   statements are equivalent:
  \begin{enumerate}
 \item The germ $(X,0)$ does not have exceptional cones.
 %\item The divisor $D= \nu_\varphi^{-1}(Z)$ of the relative Nash modification $\Na_\varphi \X  \to \X$, does not
 %      have vertical components over $Z(0)$.
 \item The pair $(\X^0, Y)$ satisfies Whitney's condition a) at the origin.
 \item The pair $(\X^0, Y)$ satisfies Whitney's conditions a) and b) at the origin.
 \item The germ $(\X,0)$ does not have exceptional cones.
\end{enumerate}
 \end{Theorem}
 \begin{proof}
     Let $(X,0)=\bigcup_{i=1}^r (X_i,0)$ be the irreducible decomposition of $(X,0)$. Then by 
  corollary \ref{CompsIrrSinConosExcepcionales}, and lemma \ref{Especializaciones} it is enough to 
  verify these equivalences for each irreducible component $(X_j,0)$ and its specialization
  space $(\X_j,0)$. Now for an irreducible germ we have:\\
  %$1) \Rightarrow 2)$ proposition \ref{VerticalComponentsCI}.\\ 
  $1) \Rightarrow 2)$ by proposition \ref{PruebamayorqueunoCI}.\\
  $2) \Rightarrow 3)$ by proposition \ref{aimplicab}.\\
  $3) \Rightarrow 4)$ by \ref{PropoGeneral*}.\\ 
  $4) \Rightarrow 1)$ by remark \ref{Remarkexcepcional}. 
 \end{proof}
  
    Suppose that  $(X,0)$ has as isolated singularity, but $C_{X,0}$ doesn't, then:
 \begin{enumerate}
  \item Either the singular locus of $\X(0)$ is contained in the singular locus of $\X$ and so this last space
 has an irreducible component contained in the special fiber $\X(0)$. 
  \item Or, every point  $p \in \mathrm{Sing} \; \X(0)$ is smooth in $\X$, which implies that the "vertical" hyperplane
  $H_t:=\{t=0\}$ is tangent to $\X$ at $p$, and so $H_t$ is a limit of tangent hyperplanes to $(\X,0)$.
 \end{enumerate}
   In any case, this will prevent us from building a Whitney stratification of $\X$ having $Y$ as a stratum.
  This kind of phenomenom is quite general and has little to do with the isolated singularity case. So
  in order to be able to build the Whitney stratification we want, it is important to have some control on
  the behavior of the singular locus of $\X$. The following lemma will help us manage this situation 
  in the case of a complete intersection tangent cone.

\begin{lema}\label{CerraduradeZ}
     Let $\nu_\varphi: \Na_\varphi(\X) \to \X $ and $(Z,0) \subset (\X,0)$ be defined by the ideal
  $I$ of \ref{BuenidealI} as before. Let $D= \nu_\varphi^{-1}(Z)$ be the exceptional divisor. If $D$ does 
  not have vertical components over $\X(0)$, then set-theoretically, the closure of $Z \setminus Z(0)$ in $\X$ 
  is equal to $Z$.
 \end{lema}
 \begin{proof}
     Let us consider the map $h:(Z,0) \to (\C,0)$ as before. If $h$ is flat, we have nothing to prove, 
  so suppose h is not flat. Then, we can find a minimal primary decomposition of 
  $I$ in $O_{\X,0}$:
  \[I = Q_1 \cap Q_2 \cap \cdots \cap Q_s\]
  such that $t^{n_i} \in Q_i$ for $1< r\leq i \leq  s$ with $n_i > 0$, so it corresponds to a
  possibly embedded irreducible component of the germ $(Z,0)$ contained in the special fiber $Z(0)$.\\
  
  	 Let $I= Q \cap B$, where $B= Q_r \cap \cdots \cap Q_s$. There exists a small 
  neighbourhood of 
  the origin $U \subset \X$, such that $I(U)= Q(U) \cap B(U)$, and for every $x \in U$
  we have the equality $I_x= Q_x \cap B_x$ in $O_{\X,x}$. But, for any open
  set $V \subset U$ such that $0 \notin V$, since $t^m \in B(V)$ and $t^m$ is a unit in $O_{\X}(V)$
  we have that $I_x=Q_x$ in $O_{\X,x}$ for any point $x \in Z\setminus \{0\}$, so
  their integral closures are equal $\overline{I_x}= \overline{Q_x}$ for every point $x \in V$\\
  
     Let $ \widetilde{\nu_\varphi}: \widetilde{\Na_\varphi(\X)}\stackrel{n}\longrightarrow  \Na_\varphi(\X) \stackrel{\nu_\varphi}
   \longrightarrow \X $ be the composition of $\nu_\varphi$ and the normalisation of $\Na_\varphi(\X)$. 
   By hypothesis, $D$ does not have
   vertical components over the origin, and since the normalisation is a finite map, we have that
   $\overline{D}= \widetilde{\nu_\varphi}^{-1}(Z)=n^{-1}(D)$ does not have vertical components over the origin either.
   Let $w \in Q$, then for $U$ sufficiently small $w \in Q(U)$. Now, we know that the coherent ideal
   $\widetilde{I}:=I O_{\widetilde{\Na_\varphi(\X)}}$ is locally invertible, so in
   particular for any point $p \in \overline{D}$ there exists an open neighborhood $V_p$ of $p$ 
   in $\widetilde{\Na_\varphi(\X)}$ such that $\widetilde{I}(V_p)=\left< g_p \right> 
   O_{\widetilde{\Na_\varphi(\X)}}(V_p)$. \\
   
     For any such neighborhood, we can consider the meromorphic function $q:=(w \circ \widetilde{\nu_\varphi} )/ g_p$. 
   The polar locus of $q$ is contained in $\overline{D}$, more precisely, since the ideal 
   $\widetilde{I}$ and $\widetilde{Q}$ coincide outside $\widetilde{\nu_\varphi}^{-1}(0)$, we have that the polar 
   locus of $q$ is contained in  $\widetilde{\nu_\varphi}^{-1}(Z(0))$. But $\overline{D}$ does not have vertical components 
   over $\X(0)$ so $\widetilde{\nu_\varphi}^{-1}(Z(0))$ is of codimension at least 2.  Since in a normal space the
   polar locus of a meromorphic function is of codimension one or empty (\cite[Thm. 71.12,p. 307]{Ka}),
   $q$ is actually holomorphic and $\widetilde{I}= \widetilde {Q}$ in $\widetilde{\Na_\varphi(\X)}$, which
   implies by \cite[Thm 2.1, p. 799]{Lej-Te} that the integral closures $\overline{I}=\overline{Q}$ 
   are equal in $O_{\X,0}$.\\
   
     Finally, since the integral closure of an ideal is contained in its radical, then
   set theoretically $Z$ is the zero locus of $\overline{I}$, that its
   $|Z| = V(\overline{I})= V(\overline{Q})=V(Q)$ and it does not have vertical components
   over the origin.\\    
  \end{proof} 

       Suppose now that both $(X,0)$ and its tangent cone are reduced complete intersections, then
   the specialization space $(\X,0)$ is also a complete intersection. In particular, refering back
   to remark \ref{ChoosingAlpha}, there is no need to choose an $\alpha$, and the ideal $I$ of 
   \ref{BuenidealI} can be chosen as the relative jacobian ideal $J_c(JM_\varphi (F))$ which
   set-theoretically defines the relative singular locus $\mathrm{Sing}_\varphi \; \X$. Note
   that that the restriction of the ideal $J_c(JM_\varphi (F))$ to the special fiber is equal to the 
   jacobian ideal $J_{C_{X,0}}$ of the tangent cone $C_{X,0}$ in $O_{C_{X,0}}$. This implies that the strict 
   transform of $\X(0)$ with respect to this blowup is equal to the Nash modification $\Na C_{X,0}$ of the fiber.

\begin{proposition}\label{Coro512}
   Let $(X,0)\subset (\C^{n+1},0)$ be a reduced germ of singularity such that the 
   tangent cone $C_{X,0}$ is a reduced complete intersection. 
   Let $|\mathrm{Sing}\, C_{X,0}|= \bigcup E_\alpha$ be the 
   irreducible decomposition of the singular locus of the tangent cone. If there exists an $\alpha$, 
   such that $E_\alpha$ is not completely contained in the reduced tangent cone $|C_{|\mathrm{
   Sing \, X|},0}|$, then it is contained in an exceptional cone. In particular we have the 
   inclusion
   \[\left|\mathrm{Sing}\, C_{X,0}\right| \subset \left|C_{|\mathrm{
   Sing \, X}|,0}\right| \bigcup \left\{ \mathrm{Exceptional \; cones}\right\}\]
 \end{proposition}
 \begin{proof}
       Let $\varphi:(\X,0) \to (\C,0)$ be the specialization space of $X$ to its tangent cone
   $(C_{X,0},0)$, and let $\nu_\varphi:\Na_\varphi(\X) \to \X $ be considered as the blowup
   of $\X$ with center $Z \subset \X$ defined by the ideal $J_c(JM_\varphi(F))$, and exceptional
   divisor $D \subset \Na_\varphi(\X) $. Since set-theoretically $Z$ is the relative singular locus
   of $\X$, then if we set $W$ as the closure of $Z \setminus Z(0)$ in $\X$, then set theoretically 
   $W(0)$ is $|C_{|\mathrm{Sing \, X}|,0}|$, so the existence of the
   $E_\alpha$ in the hypothesis amounts to $Z$ having a vertical (irreducible) component $Z_\beta$
   over the origin.\\
   
     The existence of such a $Z_\beta$ implies by \ref{CerraduradeZ} the existence of a vertical
  component $D_\beta$ of $|D|$, which then implies by \ref{VerticalComponentsCI} that the germ $(X,0)$
  has exceptional cones. Now for any point $z \in Z_\beta \setminus W$ there exists an open
  neighborhood $z \in U_z \subset \X$ such that $U_z \cap W = \emptyset$ and $Z_\beta \setminus W$ is dense
  in $Z_\beta$. That is, there exists an open neighborhood $U$ of $Z_\beta \setminus W$ in $\X$, such
  that $U \cap W = \emptyset$, and so $\nu_\varphi^{-1}(U \cap W)= \nu_\varphi^{-1}(U) \cap 
  \nu_\varphi^{-1}(W) = \emptyset$. But $\nu_\varphi^{-1}(W)$ contains $\overline{D \setminus
  D(0)}$, and $\nu_\varphi^{-1}(U) \cap D$ is not empty, so there is necessarily an irreducible 
  component $D_\beta$ of $D$, such that $D_\beta \supset \nu_\varphi^{-1}(Z_\beta)$ and $D_\beta$ 
  is completely contained in $D(0)$. All that is left to prove is that the component $D_\beta$ 
  is mapped by $\nu_\varphi$ into an exceptional cone.\\

    By remark \ref{EspacioZIntCom}, the strict transform $\overline{\nu_\varphi^{-1}(\X(0)\setminus Z)}$ is
  equal to the Nash modi- fication of the fiber $\X(0)$ which has dimension $d$, on the other hand since 
  $D_\beta$ is an irreducible component of the divisor $D$ it is also of dimension $d$ and so cannot be
  contained in $\Na \X(0)$, i.e. $D_\beta \nsubseteq \Na \X(0)$.\\

     Now, by \cite[Proposition 2.1.4.1, p. 562]{L-T2}, the cones of the aureole are set theoretically
  the images by $\kappa_\varphi$ of the irreducible components of $|\kappa_\varphi^{-1}(\X(0))|$. So
  let us consider the relative version of the diagram given in proposition \ref{ConormalvsNash},
  relating the relative Nash modification $\Na_\varphi \X$ with the relative conormal space
  $C_\varphi(\X)$.
  \[\xymatrix{
         & E_\varphi \ar @{^{(}->}[r] \ar[dr]_{p_2} \ar[dl]^{p_1} & \X \times G(n+1-d, n+1) \times \check{\P}^n   \\
     \Na_\varphi \X \ar[dr]_{\nu_\varphi} & & C_\varphi(\X) \ar[dl]^{\kappa_\varphi} \\
      & \X & }\]
  By commutativity of the diagram, we have the equality $p_2(p^{-1}(\Na \X(0)))= C(\X(0))$, where $C(\X(0))$ 
  denotes the conormal space of the fiber $\X(0)$ and it is equal to $\overline{\kappa_\varphi^{-1}(\X(0)\setminus 
  Z)}$. This implies that the space $\widetilde{D_\beta}:= p_2(p_1^{-1}(D_\beta))$ can not be contained
  in $C(\X(0))$. Now, the conormal space $C(\X(0))$ is of dimension $n$, and since $C_\varphi(\X)\to \X \to \C$ is 
  isomorphic to the specialization space of $C(X)$ to its normal cone along $\kappa_X^{-1}(0)$
  (\cite[Lemma A.4.1, p. 190]{Sab1}), then the dimension of $\kappa_\varphi^{-1}(\X(0))$ is also $n$. 
  This means that $\widetilde{D_\beta}$ is contained in an irreducible component of $|\kappa_\varphi^{-1}(\X(0))|$
  outside of $C(\X(0))$ and so is mapped by $\kappa_\varphi$ into an exceptional cone.
 \end{proof}

   Note that we always have the inclusion $\left|C_{|\mathrm{ Sing \, X}|,0}\right| \subset 
  \left|\mathrm{Sing}\, C_{X,0}\right|$, so the absence of exceptional cones together
  with \ref{Coro512} tells us that in this setting the relative singular locus, and the singular
  locus of $\X$ coincide. In particular we have $\left|C_{|\mathrm{ Sing \, X}|,0}\right| = 
  \left|\mathrm{Sing}\, C_{X,0}\right|$ and this leaves us in a good position to continue building a Whitney
  stratification of $\X$ having $Y$ as a stratum. 

 \begin{corollary}\label{AlreadyStratification} 
   Let $(X,0)$ satisfy the hypothesis of theorem \ref{Equivalenciashipersuperficies}. If $(X,0)$ has
  an isolated singularity and its tangent cone is a complete intersection singularity,
  then the absence of exceptional cones implies that $C_{X,0}$ has an isolated singularity and  
   $\{\X \setminus Y, Y \}$ is a Whitney stratification of $\X$.  
 \end{corollary}
 \begin{proof}
   Proposition \ref{Coro512} tells us that 
 $\left|C_{|\mathrm{ Sing \, X}|,0}\right| = \left|\mathrm{Sing}\, C_{X,0}\right|$, and since 
 $(X,0)$ has an isolated singularity then $\left|C_{|\mathrm{ Sing \, X}|,0}\right| = \{0\}$ and
 so the tangent cone $(C_{X,0},0)$ also has an isolated singularity. This implies, 
 that $\mathrm{Sing} \X = Y$, and theorem \ref{Equivalenciashipersuperficies} finishes the proof.
 \end{proof}

   There is a partial converse to the corollary, in which we can construct a Whitney stratification of $\X$ 
  under the assumption that the tangent cone has an isolated singularity at the origin. 

\begin{corollary}\label{Bastaconotangenteliso}
    Let $(X,0)$ satisfy the hypothesis of theorem \ref{Equivalenciashipersuperficies}. If the
 tangent cone $(C_{X,0},0)$ has an isolated singularity at the origin, then $(X,0)$ has an 
 isolated singularity and $\{\X \setminus Y, Y \}$ is a Whitney stratification of $\X$. 
\end{corollary}
\begin{proof}
   The first step is to prove that $(X,0)$ doesn't have exceptional cones, however by 
 \cite[Prop. 2.1.4.2, p. 563]{L-T2} this is always the case when the tangent cone has an isolated singularity 
 at the origin. 

   Now, by theorem \ref{Equivalenciashipersuperficies}, it is enough to prove that the singular locus
 of $\X$ is $Y$. It is a general fact that the relative singular locus  $\mathrm{Sing}_\varphi \X$ of $\X$, 
 contains the singular locus $\mathrm{Sing} \X$, and they coincide away from the special fiber. 
 In other words, the space 
 $W:= \mathrm{Sing}_\varphi \X \setminus \{\X(0)\}$ is isomorphic via $\phi: \X \setminus \X(0) 
 \to X \times \C^*$ to $\mathrm{Sing} X \times \C^*$, and so the map induced by $\varphi$ to its closure 
 $\overline{W} \to \C$ can be identified with the specialization space of $|\mathrm{Sing} X|$ to its
 tangent cone. In view of this, the hypothesis tells us that the only singular 
 point of $\X$ in the special fiber is the origin $(0,0)$; this implies $\overline{W}(0)=\{0\}$ and since 
 it is isomorphic to the tangent cone $C_{|\mathrm{Sing} X|, 0}$, then $(X,0)$ has an isolated singularity and 
 $\mathrm{Sing} \X =Y$ which finishes the proof.  
\end{proof}

\begin{example}
   Let $(V,0) \subset (\C^{n+1},0)$ be a reduced and irreducible isolated complete intersection variety defined by an homogeneous
 ideal $I_0=\left< h_{m_1}, \ldots, h_{m_k}\right>$, where $m_i$ is the degree of the polynomial. 
 That is, $V$ is the cone over a smooth, complete intersection, projective variety.\\

  Let $(X,0)\subset (\C^{n+1},0)$ be the germ defined by the ideal
 $I=\left<h_1,\ldots,h_k\right>$, where $h_i= h_{m_i}+ P_i$ and $P_i \in \C\{z_0, \ldots , z_n\}$ is such that 
 $\mathrm{ord}_0 \; P_i(z) > m_i$. Then:
 \begin{itemize}
  \item The germ $(X,0)$ is a reduced complete intersection.
  \item The tangent cone $C_{X,0}$ is defined by the ideal $I_0$ and so it is isomorphic to $V$.
 \end{itemize}

   That $X$ is a complete intersection can be seen by considering the analytic family $\{X_t\}$ defined by
  the $h_i^t:= h_{m_i} + t P_i$ and the upper semicontinuity of fiber dimension. For the other
  assertion consider the radical idea $\widetilde{I}:= \sqrt{I}$ defining $|X|$. This gives us the following 
  inclusion of initial ideals 
    \[\mathrm{In}_{\M}I_0=I_0 \subset\mathrm{In}_{\M}I \subset\mathrm{In}_{\M}\widetilde{I}\]
  and as a result the surjective morphism of analytic algebras:
   \begin{align*}
    \frac{\C\{z_0, \ldots ,z_n\}}{I_0} & \longrightarrow\frac{\C\{z_0, \ldots ,z_n\}}{\mathrm{In}_{\M}\widetilde{I}}\\ 
           O_{V,0} & \longrightarrow O_{C_{|X|,0}}
 \end{align*}
  But $V$ is irreducible, so $O_{V,0}$ is an integral domain and since both algebras have krull dimension $n+1-k$
  they are isomorphic and $I_0=\mathrm{In}_{\M}\widetilde{I}$. Finally, this tells us that 
  $\mathrm{In}_{\M}\widetilde{I}=\left< \mathrm{In}_{\M}h_1, \ldots, \mathrm{In}_{\M}h_k\right>$,
  which implies that $\widetilde{I}=\left< h_1, \ldots, h_k\right>=I$ and so $X$ is reduced and
  $C_{X,0}=V$. 
  
    Now, by construction, the specialization space $\varphi:\X \to \C$ is defined by the equations 
 $H_i(z,t)=t^{-m_i}h_i(tz)$ in $\C^{n+1} \times \C$ and since the tangent cone $C_{X,0}$ is reduced and has 
 an isolated singularity at the origin, corollary \ref{Bastaconotangenteliso} tells us that
  $\{\X \setminus Y, Y \}$ is a Whitney stratification of $\X$. 
\end{example}

\section{Conclusion}  

  We have verified that the absence of exceptional cones allows us to start building a Whitney 
 stratification of $\X$ having $Y$ as a stratum. The question now is how to continue. Proposition
 \ref{Coro512} tells us, at least in the complete intersection case, that the singular locus of $\X$ coincides
 with the specialization space $Z$ of $|\mathrm{Sing} \, X|$ to its tangent cone. 

   Suppose now, that the germ $(|\mathrm{Sing} \, X|,0)$ has a reduced tangent cone, then a stratum $\X_\lambda$ 
 containing a dense open set of $Z$ will satisfy Whitney's conditions along $Y$ if and only if the 
 germ $(|\mathrm{Sing} \, X|,0)$ doesn't have exceptional cones. 

   In view of this it seems reasonable to start by assuming the existence of a Whitney 
 stratification $\{X_\lambda\}$ of $(X,0)$ such that for every $\lambda$ the germ $(\overline{X_\lambda},0)$
 has a reduced tangent cone and no exceptional cones. In this case, the specialization space $Z_\lambda$ of 
 $\overline{X_\lambda}$ is canonically embedded as a subspace of $\X$, and the partition of $\X$ associated
 to the filtration given by the $Z_\lambda$ is a good place to start looking for the desired Whitney
 stratification of $\X$
      
\bibliographystyle{plain}

\bibliography{bibliothese}

%\begin{thebibliography}{xxxxx}

%    \bibitem[Le-Te] {Le-Te} M. Lejeune-Jalabert and B. Teissier, \textit{Cl\^oture int\'egrale
%    des ideaux et \'equisingularit\'e}, Annales de la Facult\'e des Sciences de Toulouse, Vol. XVII,
%    No. 4, 2008, 781-859. 
%  \bibitem[L-T1]{L-T1} L{\^e} D.T. and B. Teissier, \textit{Sur la g{\'e}om{\'e}trie des surfaces 
%  complexes.I.Tangentes exceptionelles}, in American Journal of Math 101, No. 2, p. 420-452.   
%  \bibitem[L-T2]{L-T2} L{\^e} D.T. and B. Teissier, \textit{Limites d'espaces tangents en
%   g{\'e}om{\'e}trie analytique}, Comm. Math. Helv., \textbf{63}, 1988, p. 540-578.
%  \bibitem[Na]{Na} V. Navarro Aznar, \textit{Conditions de Whitney et sections planes},
%   Inventiones Math.\textbf{61}, 1980, p. 199-225.   
%  \bibitem[Te1]{Te1} B. Teissier, \textit{Vari\'e\'es polaires 2: Multiplicit\'es polaires, sections 
%    planes, et conditions de Whitney}, Springer Lecture Notes \textbf{961}, 1981, p. 314-491. 
%  \bibitem[Te2]{Te2} B. Teissier, \textit{Cycles \'evanescents, secitons planes, et conditions de 
%       Whitney}, Singularit\'es \`a Carg\`ese, Ast\'erisque 7-8, 1973, 282-362.      
%  \bibitem[Whi] {Whi} H. Whitney, \textit{Tangents to an analytic variety}, Annals of Math., 
%  \textbf{81}, (1964), p. 496-549. 
%\end{thebibliography}

 \end{document}